\documentclass[letterpaper,12pt]{article}%
\usepackage{amssymb}
\usepackage{amsmath}
\usepackage{amsfonts}
\usepackage[american,english]{babel}%
\usepackage{graphicx}
\providecommand{\U}[1]{\protect\rule{.1in}{.1in}}
\begin{document}

\title{Nonparametric Estimation and Testing on Discontinuity of Positive Supported
Densities: \\A Kernel Truncation Approach\thanks{The first author gratefully acknowledges
financial support from the German Science Foundation (DFG) through the
Collaborative Research Center SFB 823. \ The second author gratefully
acknowledges financial support from Japan Society of the Promotion of Science
(grant number 15K03405).}}
\author{%
\begin{tabular}
[c]{ccc}%
Benedikt Funke\thanks{Department of Mathematics, Technical University of
Dortmund, Vogelpothsweg 87, 44227 Dortmund, Germany; e-mail:
benedikt.funke@mathematik.tu-dortmund.de.} &  & Masayuki
Hirukawa\thanks{Faculty of Economics, Setsunan University, 17-8 Ikeda
Nakamachi, Neyagawa, Osaka, 572-8508, Japan; e-mail:
hirukawa@econ.setsunan.ac.jp.}\\
Technical University &  & Setsunan University\\
of Dortmund &  &
\end{tabular}
}
\date{July 2016}
\maketitle

\begin{abstract}
\setlength{\baselineskip}{12pt}Discontinuity in density functions is of
economic importance and interest. \ For instance, in studies on regression
discontinuity designs, discontinuity in the density of a running variable
suggests violation of the no-manipulation assumption. \ In this paper we
develop estimation and testing procedures on discontinuity in densities with
positive support. \ Our approach is built on splitting the gamma kernel (Chen,
2000) into two parts at a given (dis)continuity point and constructing two
truncated kernels. \ The jump-size magnitude of the density at the point can
be estimated nonparametrically by two kernels and a multiplicative bias
correction method. \ The estimator is easy to implement, and its convergence
properties are delivered by various approximation techniques on incomplete
gamma functions. \ Based on the jump-size estimator, two versions of test
statistics for the null of continuity at a given point are also proposed.
\ Moreover, estimation theory of the entire density in the presence of a
discontinuity point is explored. \ Monte Carlo simulations confirm nice
finite-sample properties of the jump-size estimator and the test statistics. \ 

\ 

\textbf{Keywords:} boundary bias; density estimation; discontinuous
probability density; gamma kernel; incomplete gamma functions; nonparametric
kernel testing; regression discontinuity design.

\textbf{JEL Classification Codes: }C12; C13; C14.

\textbf{MSC 2010 Codes:} 62G07; 62G10; 62G20.

\end{abstract}

\begin{center}
\bigskip\setcounter{page}{0}\thispagestyle{empty}
\end{center}

\setlength{\baselineskip}{26pt}

\section{Introduction}

The objective of this paper is to develop new estimation and testing
procedures of discontinuity in density functions with support on $%
\mathbb{R}
_{+}$. \ Inference on possibly discontinuous densities has been explored in
nonparametric statistics: examples include Liebscher (1990), Cline and Hart
(1991), and Chu and Cheng (1996), to name a few. \ Discontinuity in densities
is also of economic importance and interest. \ Local randomization of a
continuous running variable is a key requirement for the validity of
regression discontinuity designs (\textquotedblleft RDD\textquotedblright); if
the value of the running variable falls into the left and right of the cutoff
strategically, then treatment effects are no longer point identified due to
self-selection. \ Therefore, detection of discontinuity in the density of the
running variable at the cutoff suggests evidence of such strategic behavior or
manipulation in RDD. \ Nonetheless, estimation and inference on jump-size
magnitudes of densities at discontinuity points have not attracted interest in
econometrics up until recently. \ McCrary (2008) applies a bin-based local
linear regression method to estimate jump sizes. \ Subsequently, Otsu, Xu and
Matsushita (2013) propose two versions of empirical likelihood-based inference
procedures grounded on binning and local likelihood methods. \ While our
proposal can be viewed as an extension of these articles, it has a unique
feature. \ In our approach, jump sizes are estimated by means of density
estimation techniques using the kernels obtained through truncating asymmetric
kernels at a given (dis)continuity point, unlike nonparametric regression or
local likelihood approaches using standard symmetric kernels.

Before proceeding, it is worth explaining why we specialize in asymmetric
kernel smoothing. \ Empirical studies on discontinuity in densities frequently
pay attention to the distributions of economic variables such as (taxable or
relative) incomes (Saez, 2010; Bertrand, Kamenica and Pan, 2015), wages
(DiNardo, Fortin and Lemieux, 1996), school enrollment counts (Angrist and
Lavy, 1999) and proportion of votes for proposed bills (McCrary, 2008). \ The
distributions, if they are free of discontinuity points, can be empirically
characterized by two stylized facts, namely, (i) existence of a lower bound in
support (most possibly at the origin) and (ii) concentration of observations
near the boundary and a long tail with sparse data. \ When estimating such
densities nonparametrically using symmetric kernels, we must rely either on a
boundary correction method and an adaptive smoothing technique (e.g., variable
bandwidth methods) \emph{simultaneously}, or on back-transforming the density
estimator from the log-transformed data to the original scale. \ The former is
apparently cumbersome, and density estimates by the latter often behave poorly
(e.g., Cowell, Ferreira and Litchfield, 1998) although the method is popularly
applied in empirical works. \ Asymmetric kernels with support on $%
\mathbb{R}
_{+}$ have emerged as a viable alternative that can accommodate the stylized
facts. \ Although there are various classes of asymmetric kernels, for the
sake of simplicity and due to popularity this study focuses exclusively on the
gamma kernel by Chen (2000) \
\[
K_{G\left(  x,b\right)  }\left(  u\right)  =\frac{u^{x/b}\exp\left(
-u/b\right)  }{b^{x/b+1}\Gamma\left(  x/b+1\right)  }\mathbf{1}\left(
u\geq0\right)  ,
\]
where $x\left(  \geq0\right)  $ and $b\left(  >0\right)  $\ are the design
point and smoothing parameter, respectively. \ 

When the density has a discontinuity point, the jump-size magnitude at the
point can be defined as the difference between left and right limits of the
density at the point. \ While nonparametric regression (McCrary, 2008) and
empirical likelihood (Otsu, Xu and Matsushita, 2013) methods have been applied
to estimate the jump size, we attempt to have our jump-size estimator preserve
appealing properties of the gamma kernel. \ Accordingly, we split the gamma
kernel into two parts at the discontinuity point, and make each part a
legitimate kernel by re-normalization. \ The left and right limits of the
density can be estimated by two truncated kernels. \ Although the estimators
are consistent and their variance convergences are usual $O\left(
n^{-1}b^{-1/2}\right)  $ where $n$ is the sample size, their bias convergences
are $O\left(  b^{1/2}\right)  $, not the usual $O\left(  b\right)  $ rate.
\ Then, we apply the multiplicative bias correction technique by Terrell and
Scott (1980) to eliminate the undesirable $O\left(  b^{1/2}\right)  $ biases
without inflating the order of magnitude in variance. \ Moreover, we take
particular care of choosing the smoothing parameter. \ Specifically, the
method of power-optimality smoothing parameter selection by Kulasekera and
Wang (1998) is tailored to inference problems on discontinuous densities.

Our proposal has three contributions to the literature. \ First, unlike the
binned local linear (\textquotedblleft BLL\textquotedblright) estimation by
McCrary (2008), our kernel truncation approach always generates nonnegative
density estimates and is free from choosing bin widths. \ Our jump-size
estimator is also easy to implement.\ Since it has a closed form, nonlinear
optimization as in Otsu, Xu and Matsushita (2013) is unnecessary. \ While
incomplete gamma functions are key ingredients in our estimator, standard
statistical packages including GAUSS, Matlab and R prepare a command that can
return values of the functions either directly or in the form of gamma
cumulative distribution functions\textbf{. }

Second, in delivering convergence results of asymmetric kernel estimators, we
utilize the mathematical tools and proof strategies that are totally different
from those for nonparametric estimators smoothed by symmetric kernels.
\ Asymptotic results throughout this paper are built upon a few different
approximation techniques on incomplete gamma functions; such proof strategies
are taken for the first time in the econometric literature, to the best of our knowledge.

Third, we also present estimation theory of the entire density in the presence
of a discontinuity point. \ Indeed, Imbens and Lemieux (2008) argue importance
of graphical analyses in empirical studies on RDD, including inspections of
densities of running variables. \ It is demonstrated that density estimators
smoothed by the truncated gamma kernels admit the same bias and variance
approximations as the gamma kernel density estimator does. \ Furthermore, the
truncated gamma-kernel density estimator is shown to be consistent even when
the true density is unbounded at the origin.

The remainder of this paper is organized as follows. \ Section 2 presents
estimation and testing procedures of the density at a known discontinuity
point $c\left(  >0\right)  $. \ As an important practical problem, a smoothing
parameter selection method is also developed. \ Our particular focus is on the
choice method for power optimality. \ In Section 3, we discuss how to estimate
the entire density when the density has a discontinuity point. \ Convergence
properties of density estimates are also explored. \ Section 4 conducts Monte
Carlo simulations to evaluate finite-sample properties of the proposed
jump-size estimator and test statistic. \ An empirical application on the
validity of RDD is presented in Section 5. \ Section 6 summarizes the main
results of the paper. \ Proofs are provided in the Appendix.

This paper adopts the following notational conventions: for $a>0$,
$\Gamma\left(  a\right)  =\int_{0}^{\infty}t^{a-1}\exp\left(  -t\right)  dt$
is the gamma function; for $a,z>0$, $\gamma\left(  a,z\right)  =\int_{0}%
^{z}t^{a-1}\exp\left(  -t\right)  dt$\ and $\Gamma\left(  a,z\right)
=\int_{z}^{\infty}t^{a-1}\exp\left(  -t\right)  dt=\Gamma\left(  a\right)
-\gamma\left(  a,z\right)  $ denote the lower and upper incomplete gamma
functions, respectively; $\mathbf{1}\left\{  \cdot\right\}  $ signifies an
indicator function; and $\left\lfloor \cdot\right\rfloor $ denotes the integer
part.\ \ Lastly, the expression `$X_{n}\sim Y_{n}$' is used whenever
$X_{n}/Y_{n}\rightarrow1$\ as $n\rightarrow\infty$.

\section{Estimation and Inference for Discontinuity in the Density}

\subsection{Setup}

Suppose that we suspect discontinuity of the probability density function
(\textquotedblleft pdf\textquotedblright) $f\left(  x\right)  $ at a given
point $x=c\left(  >0\right)  $, which is assumed to be interior throughout.
\ Also let \
\[
f_{-}\left(  c\right)  :=\lim_{x\uparrow c}f\left(  x\right)  \text{ and
}f_{+}\left(  c\right)  :=\lim_{x\downarrow c}f\left(  x\right)  ,
\]
be the lower and upper limits of the pdf at $x=c$, respectively. \ Our
parameter of interest is the jump-size magnitude of the density at $c$%
\[
J\left(  c\right)  :=f_{+}\left(  c\right)  -f_{-}\left(  c\right)  .
\]
To check whether $f$ is (dis)continuous at $c$, we first estimate $J\left(
c\right)  $ nonparametrically and then proceed to a hypothesis testing for the
null of continuity of $f$ at $c$, i.e., $H_{0}:J\left(  c\right)  =0$, against
the two-sided alternative.

\subsection{An Issue in Estimating Two Limits of the Density}

To develop a consistent estimator of $J\left(  c\right)  $, we start our
analysis from estimating two limits of the density at $c$. \ Let $\left\{
X_{i}\right\}  _{i=1}^{n}$\ be a univariate random sample drawn from a
distribution that has the pdf $f$. \ When $f$ is indeed discontinuous at $c$,
a reasonable method would be to estimate $f_{-}\left(  c\right)  $ and
$f_{+}\left(  c\right)  $ using sub-samples $\left\{  X_{i}^{-}\right\}
:=\left\{  X_{i}:X_{i}<c\right\}  $ and $\left\{  X_{i}^{+}\right\}
:=\left\{  X_{i}:X_{i}\geq c\right\}  $, respectively. \ Instead of relying on
nonparametric regression or local likelihood methods, we split the gamma
kernel into two parts at $c$, namely,
\[
K_{G\left(  x,b\right)  }\left(  u\right)  :=K_{G\left(  x,b;c\right)  }%
^{L}\left(  u\right)  +K_{G\left(  x,b;c\right)  }^{U}\left(  u\right)  ,
\]
where
\begin{align*}
K_{G\left(  x,b;c\right)  }^{L}\left(  u\right)   & =\frac{u^{x/b}\exp\left(
-u/b\right)  }{b^{x/b+1}\Gamma\left(  x/b+1\right)  }\mathbf{1}\left(  0\leq
u<c\right)  \text{ and}\\
K_{G\left(  x,b;c\right)  }^{U}\left(  u\right)   & =\frac{u^{x/b}\exp\left(
-u/b\right)  }{b^{x/b+1}\Gamma\left(  x/b+1\right)  }\mathbf{1}\left(  u\geq
c\right)  .
\end{align*}
However, neither $K_{G\left(  x,b;c\right)  }^{L}\left(  u\right)  $ nor
$K_{G\left(  x,b;c\right)  }^{U}\left(  u\right)  $\ is a legitimate kernel
function in the sense that%
\begin{align*}
\int_{0}^{\infty}K_{G\left(  x,b;c\right)  }^{L}\left(  u\right)  du  &
=\frac{\gamma\left(  x/b+1,c/b\right)  }{\Gamma\left(  x/b+1\right)  }\text{
and}\\
\int_{0}^{\infty}K_{G\left(  x,b;c\right)  }^{U}\left(  u\right)  du  &
=\frac{\Gamma\left(  x/b+1,c/b\right)  }{\Gamma\left(  x/b+1\right)  }.
\end{align*}
Therefore, we make scale-adjustments to obtain the re-normalized truncated
kernels as
\begin{align*}
K_{G\left(  x,b;c\right)  }^{-}\left(  u\right)   & =\frac{\Gamma\left(
x/b+1\right)  }{\gamma\left(  x/b+1,c/b\right)  }K_{G\left(  x,b;c\right)
}^{L}\left(  u\right)  \text{ and }\\
K_{G\left(  x,b;c\right)  }^{+}\left(  u\right)   & =\frac{\Gamma\left(
x/b+1\right)  }{\Gamma\left(  x/b+1,c/b\right)  }K_{G\left(  x,b;c\right)
}^{U}\left(  u\right)  .
\end{align*}
These kernels yield estimators of $f_{-}\left(  c\right)  $ and $f_{+}\left(
c\right)  $\ as%
\begin{align*}
\hat{f}_{-}\left(  c\right)   & =\frac{1}{n}\sum_{i=1}^{n}\left.  K_{G\left(
x,b;c\right)  }^{-}\left(  X_{i}\right)  \right\vert _{x=c}=\frac{1}{n}%
\sum_{i=1}^{n}K_{G\left(  c,b;c\right)  }^{-}\left(  X_{i}\right)  \text{ and
}\\
\hat{f}_{+}\left(  c\right)   & =\frac{1}{n}\sum_{i=1}^{n}\left.  K_{G\left(
x,b;c\right)  }^{+}\left(  X_{i}\right)  \right\vert _{x=c}=\frac{1}{n}%
\sum_{i=1}^{n}K_{G\left(  c,b;c\right)  }^{+}\left(  X_{i}\right)  .
\end{align*}

To explore asymptotic properties of these estimators, we make the following
assumptions. \ For notational conciseness, expressions such as
\textquotedblleft$f_{\pm}\left(  c\right)  $\textquotedblright\ are used
throughout, whenever no confusions may occur.

\begin{description}
\item[Assumption 1.] The random sample $\left\{  X_{i}\right\}  _{i=1}^{n}$ is
drawn from a univariate distribution with a pdf $f$\ having support on $%
\mathbb{R}
_{+}$.

\item[Assumption 2.] The second-order derivative of the pdf $f$ is
H\"{o}lder-continuous of order $\varsigma\in\left(  0,1\right]  $ on $%
\mathbb{R}
_{+}\backslash\left\{  c\right\}  $. \ Also let $f_{-}^{\left(  j\right)
}\left(  c\right)  :=\lim_{x\uparrow c}d^{j}f\left(  x\right)  /dx^{j}$ and
$f_{+}^{\left(  j\right)  }\left(  c\right)  :=\lim_{x\downarrow c}%
d^{j}f\left(  x\right)  /dx^{j}$\ for $j=1,2$. \ Then, $f_{\pm}\left(
c\right)  >0$ and $\left\vert f_{\pm}^{\left(  2\right)  }\left(  c\right)
\right\vert <\infty$.

\item[Assumption 3.] The smoothing parameter $b\left(  =b_{n}>0\right)  $
satisfies $b+\left(  nb\right)  ^{-1}\rightarrow0$ as $n\rightarrow\infty$.
\end{description}

Assumptions 1 and 3 are standard in the literature on asymmetric kernel
smoothing (e.g., Chen, 2000; Hirukawa and Sakudo, 2015). \ The condition
\textquotedblleft$\left(  nb\right)  ^{-1}\rightarrow0$\textquotedblright\ in
Assumption 3 is required for the estimation of the entire density that will be
discussed in Section 3, whereas a weaker condition \textquotedblleft$\left(
nb^{1/2}\right)  ^{-1}\rightarrow0$\textquotedblright\ suffices for
Propositions 1 and 2 and Theorem 1 below. \ Moreover, an equivalent to
Assumption 2 can be found in McCrary (2008) and Otsu, Xu and Matsushita
(2013). \ In particular, H\"{o}lder-continuity of the second-order density
derivative $f^{\left(  2\right)  }\left(  \cdot\right)  $ in Assumption 2
implies that there is a constant $L\in\left(  0,\infty\right)  $ such that%
\begin{align*}
\left\vert f^{\left(  2\right)  }\left(  s\right)  -f^{\left(  2\right)
}\left(  t\right)  \right\vert  & \leq L\left\vert s-t\right\vert ^{\varsigma
},\,\forall s,t\in\left[  0,c\right)  \text{ and }\\
\left\vert f^{\left(  2\right)  }\left(  s^{\prime}\right)  -f^{\left(
2\right)  }\left(  t^{\prime}\right)  \right\vert  & \leq L\left\vert
s^{\prime}-t^{\prime}\right\vert ^{\varsigma},\,\forall s^{\prime},t^{\prime
}\in\left[  c,\infty\right)  .
\end{align*}

The proposition below refers to bias and variance approximations of $\hat
{f}_{\pm}\left(  c\right)  $. \ It is worth emphasizing that all convergences
results in this paper are built upon a few different approximation techniques
on incomplete gamma functions; such proof strategies are taken for the first
time in the econometric literature, to the best of our knowledge. \ Moreover,
for the purpose of our subsequent analysis, the bias expansion is derived up
to the second-order term. \ 

\paragraph{\emph{Proposition 1.}}

\textit{Under Assumptions 1-3, as }$n\rightarrow\infty$\textit{, }%
\begin{align*}
Bias\left\{  \hat{f}_{\pm}\left(  c\right)  \right\}   & \sim\mp\sqrt{\frac
{2}{\pi}}c^{1/2}f_{\pm}^{\left(  1\right)  }\left(  c\right)  b^{1/2}+\left\{
\left(  1-\frac{4}{3\pi}\right)  f_{\pm}^{\left(  1\right)  }\left(  c\right)
+\frac{c}{2}f_{\pm}^{\left(  2\right)  }\left(  c\right)  \right\}
b,\text{\textit{\ and}}\\
Var\left\{  \hat{f}_{\pm}\left(  c\right)  \right\}   & \sim\frac{1}{nb^{1/2}%
}\frac{f_{\pm}\left(  c\right)  }{\sqrt{\pi}c^{1/2}}.
\end{align*}

Proposition 1 implies that $\hat{f}_{\pm}\left(  c\right)  $ are consistent
for $f_{\pm}\left(  c\right)  $, and that their variance convergence has a
usual rate of $O\left(  n^{-1}b^{-1/2}\right)  $. \ Nevertheless, the bias
convergence is $O\left(  b^{1/2}\right)  $, which is slower than the usual
$O\left(  b\right)  $ rate. \ This is an outcome of one-sided smoothing. \ If
$f $ were continuous at $c$ and smoothing were made on both sides of the
design point $c$ using the gamma kernel, the nearly symmetric shape of the
kernel would cancel out the $O\left(  b^{1/2}\right)  $ bias.\footnote{This
can be also seen by combining two estimators $\hat{f}_{\pm}\left(  c\right)  $
as a weighted sum.} \ In reality, because data points used for estimating
$f_{\pm}\left(  c\right)  $ lie only on either the left or right side of $c$,
the $O\left(  b^{1/2}\right)  $ bias never vanishes. \ It follows that when
$J\left(  c\right)  $\ is estimated by $\hat{J}\left(  c\right)  :=\hat{f}%
_{+}\left(  c\right)  -\hat{f}_{-}\left(  c\right)  $, it also has an inferior
$O\left(  b^{1/2}\right)  $ bias. \ Therefore, our goal is to propose an
estimator of $J\left(  c\right)  $\ with an $O\left(  b\right)  $ bias and an
$O\left(  n^{-1}b^{-1/2}\right)  $ variance.

\subsection{Bias-Corrected Estimation and Inference \ }

To improve the bias convergence in estimators of $f_{\pm}\left(  c\right)
$\ from $O\left(  b^{1/2}\right)  $ to $O\left(  b\right)  $ while the order
of magnitude in variance remains unchanged, we propose to employ a
multiplicative bias correction (\textquotedblleft MBC\textquotedblright)
technique. \ As in Hirukawa (2010), Hirukawa and Sakudo (2014, 2015), and
Funke and Kawka (2015), the MBC method proposed by Terrell and Scott (1980) is
adopted.\footnote{Aforementioned articles also apply another MBC method
proposed by Jones, Linton and Nielsen (1995). \ However, it appears that the
method fails to eliminate the $O\left(  b^{1/2}\right)  $ bias. \ Their MBC
estimator of $f_{-}\left(  c\right)  $, for example, can be written as
\[
\breve{f}_{-}\left(  c\right)  :=\hat{f}_{-}\left(  c\right)  \breve{\alpha
}_{-}\left(  c\right)  :=\hat{f}_{-}\left(  c\right)  \left\{  \frac{1}{n}%
\sum_{i=1}^{n}\frac{K_{G\left(  c,b;c\right)  }^{-}\left(  X_{i}\right)
}{\hat{f}_{-}\left(  X_{i}\right)  }\right\}  ,
\]
where $\breve{\alpha}_{-}\left(  c\right)  $ serves as the `bias correction'
term. \ However, $\hat{f}_{-}\left(  x\right)  \,\left(  x<c\right)  $ has an
$O\left(  b\right)  $ bias, as stated in Theorem 2, so does $\breve{\alpha
}_{-}\left(  c\right)  $. \ Therefore, the $O\left(  b^{1/2}\right)  $ bias in
$\hat{f}_{-}\left(  c\right)  $ never vanishes, and thus we do not pursue this
type of MBC.} \ The method eliminates the leading bias term by constructing a
multiplicative combination of two density estimators with different smoothing
parameters. \ In our context, for some constant $\delta\in\left(  0,1\right)
$, the MBC estimators of $f_{\pm}\left(  c\right)  $ can be defined as%
\[
\tilde{f}_{\pm}\left(  c\right)  =\left\{  \hat{f}_{\pm,b}\left(  c\right)
\right\}  ^{1/\left(  1-\delta^{1/2}\right)  }\left\{  \hat{f}_{\pm,b/\delta
}\left(  c\right)  \right\}  ^{-\delta^{1/2}/\left(  1-\delta^{1/2}\right)  },
\]
where $\hat{f}_{\bullet,b}\left(  x\right)  $ and $\hat{f}_{\bullet,b/\delta
}\left(  x\right)  $ signify the density estimators using smoothing parameters
$b$ and $b/\delta$, respectively. \ Not only are $\tilde{f}_{\pm}\left(
c\right)  $ nonnegative by construction, but also their bias and variance
convergences are usual $O\left(  b\right)  $ and $O\left(  n^{-1}%
b^{-1/2}\right)  $ rates, respectively, as documented in the next proposition.
\ The proof is similar to the one for Theorem 1 of Hirukawa and Sakudo (2014),
and thus it is omitted.

\paragraph{\emph{Proposition 2.}}

\textit{Under Assumptions 1-3, as }$n\rightarrow\infty$,\textit{\ }%
\begin{align*}
Bias\left\{  \tilde{f}_{\pm}\left(  c\right)  \right\}   & \sim\left(
\frac{1}{\delta^{1/2}}\right)  \left[  \frac{c}{\pi}\left\{  \frac{\left(
f_{\pm}^{\left(  1\right)  }\left(  c\right)  \right)  ^{2}}{f_{\pm}\left(
c\right)  }\right\}  -\left\{  \left(  1-\frac{4}{3\pi}\right)  f_{\pm
}^{\left(  1\right)  }\left(  c\right)  +\frac{c}{2}f_{\pm}^{\left(  2\right)
}\left(  c\right)  \right\}  \right]  b,\text{\textit{\ and}}\\
Var\left\{  \tilde{f}_{\pm}\left(  c\right)  \right\}   & \sim\frac
{1}{nb^{1/2}}\lambda\left(  \delta\right)  \frac{f_{\pm}\left(  c\right)
}{\sqrt{\pi}c^{1/2}},
\end{align*}
\textit{where}
\[
\lambda\left(  \delta\right)  :=\frac{\left(  1+\delta^{3/2}\right)  \left(
1+\delta\right)  ^{1/2}-2\sqrt{2}\delta}{\left(  1+\delta\right)
^{1/2}\left(  1-\delta^{1/2}\right)  ^{2}}%
\]
\textit{is monotonously increasing in }$\delta\in\left(  0,1\right)
$\textit{\ with }%
\[
\lim_{\delta\downarrow0}\lambda\left(  \delta\right)  =1\text{\textit{\ and }%
}\lim_{\delta\uparrow1}\lambda\left(  \delta\right)  =\frac{11}{4}.
\]

Proposition 2 suggests that as $\delta\downarrow0$ ($\delta\uparrow1$) or in
case of oversmoothing (undersmoothing), the bias increases (decreases) and the
variance decreases (increases). \ It is a common practice in nonparametric
kernel testing that the bias is made asymptotically negligible via
undersmoothing, and thus what matters for inference is the size of
$\lambda\left(  \delta\right)  $. \ Because of no minimum in $\lambda\left(
\delta\right)  $, the choice of $\delta$\ is left as an exercise in Monte
Carlo simulations. \ 

It also follows that $J\left(  c\right)  $ can be consistently estimated as
$\tilde{J}\left(  c\right)  :=\tilde{f}_{+}\left(  c\right)  -\tilde{f}%
_{-}\left(  c\right)  $. \ The next theorem refers to the limiting
distribution of $\tilde{J}\left(  c\right)  $.\footnote{It is possible to use
different constants $\delta_{-}$ and $\delta_{+}$ and/or different smoothing
parameters $b_{-}$ and $b_{+}$ for $\tilde{f}_{-}\left(  c\right)  $ and
$\tilde{f}_{+}\left(  c\right)  $, as long as $b_{-}$ and $b_{+}$ shrink to
zero at the same rate. \ For convenience, however, we choose to employ the
same $\delta$ and $b$.}

\paragraph{\emph{Theorem 1.}}

\textit{Under Assumptions 1-3, as }$n\rightarrow\infty$,\textit{\ }%
\begin{equation}
\sqrt{nb^{1/2}}\left\{  \tilde{J}\left(  c\right)  -J\left(  c\right)
-B\left(  c\right)  b+o\left(  b\right)  \right\}  \overset{d}{\rightarrow
}N\left(  0,V\left(  c\right)  \right)  ,\label{Thm1}%
\end{equation}
\textit{where }%
\begin{align*}
B\left(  c\right)   & =\left(  \frac{1}{\delta^{1/2}}\right)  \left[  \frac
{c}{\pi}\left\{  \frac{\left(  f_{+}^{\left(  1\right)  }\left(  c\right)
\right)  ^{2}}{f_{+}\left(  c\right)  }-\frac{\left(  f_{-}^{\left(  1\right)
}\left(  c\right)  \right)  ^{2}}{f_{-}\left(  c\right)  }\right\}  \right. \\
& \left.  -\left\{  \left(  1-\frac{4}{3\pi}\right)  \left(  f_{+}^{\left(
1\right)  }\left(  c\right)  -f_{-}^{\left(  1\right)  }\left(  c\right)
\right)  +\frac{c}{2}\left(  f_{+}^{\left(  2\right)  }\left(  c\right)
-f_{-}^{\left(  2\right)  }\left(  c\right)  \right)  \right\}  \right] \\
V\left(  c\right)   & =\lambda\left(  \delta\right)  \left\{  \frac
{f_{+}\left(  c\right)  +f_{-}\left(  c\right)  }{\sqrt{\pi}c^{1/2}}\right\}
,
\end{align*}
\textit{and }$\lambda\left(  \delta\right)  $\textit{\ is defined in
Proposition 2. \ In addition, if }$nb^{5/2}\rightarrow0$\textit{\ as
}$n\rightarrow\infty$\textit{, then (\ref{Thm1}) reduces to}%
\[
\sqrt{nb^{1/2}}\left\{  \tilde{J}\left(  c\right)  -J\left(  c\right)
\right\}  \overset{d}{\rightarrow}N\left(  0,V\left(  c\right)  \right)  .
\]

As indicated in Proposition 2, $\tilde{J}\left(  c\right)  $ has an $O\left(
b\right)  $ bias and an $O\left(  n^{-1}b^{-1/2}\right)  $ variance. \ Observe
that for a given $\delta$, the variance coefficient decreases as $c$
increases, i.e., as the discontinuity point moves away from the origin. \ We
can also find that the leading bias term $B\left(  c\right)  b$\ cancels out
if $f$ has a continuous second-order derivative at $c$. \ 

Theorem 1 also implies that given a smoothing parameter $b=Bn^{-q}$ for some
constants $B\in\left(  0,\infty\right)  $ and $q\in\left(  2/5,1\right)  $ and
$\tilde{V}\left(  c\right)  $,\ a consistent estimate of $V\left(  c\right)
$, the test statistic is%
\[
T\left(  c\right)  :=\frac{\sqrt{nb^{1/2}}\tilde{J}\left(  c\right)  }%
{\sqrt{\tilde{V}\left(  c\right)  }}\overset{d}{\rightarrow}N\left(
0,1\right)  \text{ under }H_{0}:J\left(  c\right)  =0.
\]
Moreover, as documented in the next proposition, the test is consistent.
\ Observe that the power approaches one for local alternatives with
convergence rates no faster than $n^{1/2}b^{1/4}$, as well as for fixed alternatives.

\paragraph{\emph{Proposition 3.}}

\textit{Under Assumptions 1-3, as }$n\rightarrow\infty$, $\Pr\left\{
\left\vert T\left(  c\right)  \right\vert >B_{n}\right\}  \rightarrow
1$\textit{\ under }$H_{1}:J\left(  c\right)  \neq0$ \textit{for any
non-stochastic sequence }$B_{n}$\textit{\ satisfying }$B_{n}=o\left(
n^{1/2}b^{1/4}\right)  $\textit{.}

Our remaining tasks are to present examples of $\tilde{V}\left(  c\right)  $
and to propose a choice method of $b$. \ The latter is discussed in the next
section, whereas there are a few candidates of $\tilde{V}\left(  c\right)  $.
\ Replacing $f_{\pm}\left(  c\right)  $ in $V\left(  c\right)  $\ with their
consistent estimates $\tilde{f}_{\pm}\left(  c\right)  $\ immediately yields%
\[
\tilde{V}_{1}\left(  c\right)  :=\lambda\left(  \delta\right)  \left\{
\frac{\tilde{f}_{+}\left(  c\right)  +\tilde{f}_{-}\left(  c\right)  }%
{\sqrt{\pi}c^{1/2}}\right\}  .
\]
Alternatively, it is possible to compute the gamma kernel density estimator at
$c$
\[
\hat{f}\left(  c\right)  :=\frac{1}{n}\sum_{i=1}^{n}\left.  K_{G\left(
x,b\right)  }\left(  X_{i}\right)  \right\vert _{x=c}=\frac{1}{n}\sum
_{i=1}^{n}K_{G\left(  c,b\right)  }\left(  X_{i}\right)  .
\]
By (\ref{UF4}) and (\ref{P1-2}), we have
\[
\frac{\gamma\left(  c/b+1,c/b\right)  }{\Gamma\left(  c/b+1\right)  }=\frac
{1}{2}+O\left(  b^{1/2}\right)  \text{ and }\frac{\Gamma\left(
c/b+1,c/b\right)  }{\Gamma\left(  c/b+1\right)  }=\frac{1}{2}+O\left(
b^{1/2}\right)  .
\]
It follows that
\[
\hat{f}\left(  c\right)  =\frac{\gamma\left(  c/b+1,c/b\right)  }%
{\Gamma\left(  c/b+1\right)  }\hat{f}_{-}\left(  c\right)  +\frac
{\Gamma\left(  c/b+1,c/b\right)  }{\Gamma\left(  c/b+1\right)  }\hat{f}%
_{+}\left(  c\right)  \overset{p}{\rightarrow}\frac{f_{+}\left(  c\right)
+f_{-}\left(  c\right)  }{2}.
\]
As a consequence, we can obtain another estimator of $V\left(  c\right)  $ as
\[
\tilde{V}_{2}\left(  c\right)  :=\lambda\left(  \delta\right)  \left\{
\frac{2\hat{f}\left(  c\right)  }{\sqrt{\pi}c^{1/2}}\right\}  .
\]
\ 

\subsection{Smoothing Parameter Selection}

How to choose the value of the smoothing parameter $b$ is an important
practical problem. \ McCrary (2008) proposes the choice method which closely
follows the literature on the BLL smoothing. \ Moreover, in the literature on
RDD, Imbens and Kalyanaraman (2012) and Porter and Yu (2015, Section 5.4)
discuss methods of choosing the smoothing parameter. \ All these proposals
rely on either a cross-validation criterion or a plug-in approach, and thus
they stand on the idea of estimation-optimality. \ However, once our priority
is given to testing for continuity of the pdf $f$ at a given point $c$, such
approaches cannot be justified in theory or practice, because
estimation-optimal values may not be equally optimal for testing purposes.
\ Here we have a preference for test-optimality and thus adopt the
power-optimality criterion by Kulasekera and Wang (1998), whose idea is also
applied in Hirukawa and Sakudo (2016). \ 

Below Procedure 1 of Kulasekera and Wang (1998) is tailored to our context.
\ The procedure is a version of sub-sampling. \ Let $n_{-}$ and $n_{+}$ be the
numbers of observations in sub-samples $\left\{  X_{i}^{-}\right\}  $ and
$\left\{  X_{i}^{+}\right\}  $, respectively, where $n\equiv n_{-}+n_{+}$.
\ Also assume that $\left\{  X_{i}^{-}\right\}  _{i=1}^{n_{-}}$ and $\left\{
X_{i}^{+}\right\}  _{i=1}^{n_{+}}$ are ordered samples. \ Then, the entire
sample $\left\{  X_{i}\right\}  _{i=1}^{n}=\left\{  \left\{  X_{i}%
^{-}\right\}  _{i=1}^{n_{-}},\left\{  X_{i}^{+}\right\}  _{i=1}^{n_{+}%
}\right\}  $ can be split into $M$ sub-samples, where $M=M_{n}$ is a
non-stochastic sequence that satisfies $1/M+M/n\rightarrow0$ as $n\rightarrow
\infty$. \ Given such $M$, $\left(  k_{-},k_{+}\right)  :=\left(  \left\lfloor
n_{-}/M\right\rfloor ,\left\lfloor n_{+}/M\right\rfloor \right)  $ and
$k:=k_{-}+k_{+}$, the $m$th sub-sample is defined as\
\[
\left\{  X_{m,i}\right\}  _{i=1}^{k}:=\left\{  \left\{  X_{m+\left(
i-1\right)  M}^{-}\right\}  _{i=1}^{k_{-}},\left\{  X_{m+\left(  i-1\right)
M}^{+}\right\}  _{i=1}^{k_{+}}\right\}  ,\,m=1,\ldots,M.\
\]

The test statistic using the $m$th sub-sample $\left\{  X_{m,i}\right\}
_{i=1}^{k}$ becomes
\[
T_{m}\left(  c\right)  :=\frac{\sqrt{kb^{1/2}}\tilde{J}_{m}\left(  c\right)
}{\sqrt{\tilde{V}_{m}\left(  c\right)  }},\,m=1,\ldots,M,
\]
where $\tilde{J}_{m}\left(  c\right)  $ and $\tilde{V}_{m}\left(  c\right)  $
(which is either $\tilde{V}_{1,m}\left(  c\right)  $ or $\tilde{V}%
_{2,m}\left(  c\right)  $) are the sub-sample analogues of $\tilde{J}\left(
c\right)  $ and $\tilde{V}\left(  c\right)  $, respectively. \ Also denote the
set of admissible values for $b=b_{n}$ as $H_{n}:=\left[  \underline{B}%
n^{-q},\overline{B}n^{-q}\right]  $ for some prespecified exponent
$q\in\left(  2/5,1\right)  $ and two constants $0<\underline{B}<\overline
{B}<\infty$. \ Moreover, let
\[
\hat{\pi}_{M}\left(  b_{k}\right)  :=\frac{1}{M}\sum_{m=1}^{M}\mathbf{1}%
\left\{  T_{m}\left(  c\right)  >c_{m}\left(  \alpha\right)  \right\}  ,
\]
where $c_{m}\left(  \alpha\right)  $ is the critical value for the size
$\alpha$\ test using the $m$th sub-sample. \ We pick the power-maximized
$\hat{b}_{k}=\hat{B}k^{-q}=\arg\max_{b_{k}\in H_{k}}\hat{\pi}_{M}\left(
b_{k}\right)  $, and the smoothing parameter value $\hat{b}_{n}:=\hat{B}n^{-q}
$ follows.

We conclude this section by stating how to obtain $\hat{b}_{n}$\ in practice.
\ Step 1 reflects that $M$ should be divergent but smaller than both $n_{-}$
and $n_{+}$ in finite samples. \ Step 3 follows from the implementation
methods in Kulasekera and Wang (1998). \ Finally, Step 4 corresponds to the
case for more than one maximizer of $\hat{\pi}_{M}\left(  b_{k}\right)  $.

\begin{center}
\textbf{\ }%
\begin{tabular}
[c]{ll}%
\textbf{Step 1:} & Choose some $p\in\left(  0,1\right)  $ and specify
$M=\left\lfloor \min\left\{  n_{-}^{p},n_{+}^{p}\right\}  \right\rfloor $.\\
\textbf{Step 2:} & Make $M$ sub-samples of sizes $\left(  k_{-},k_{+}\right)
=\left(  \left\lfloor n_{-}/M\right\rfloor ,\left\lfloor n_{+}/M\right\rfloor
\right)  $.\\
\textbf{Step 3:} & Pick two constants $0<\underline{H}<\overline{H}<1$ and
define $H_{k}=\left[  \underline{H},\overline{H}\right]  $.\\
\textbf{Step 4:} & Set $c_{m}\left(  \alpha\right)  \equiv z_{\alpha}$ and
find $\hat{b}_{k}=\inf\left\{  \arg\max_{b_{k}\in H_{k}}\hat{\pi}_{M}\left(
b_{k}\right)  \right\}  $\\
& by a grid search.\\
\textbf{Step 5:} & Recover $\hat{B}$\ by $\hat{B}=\hat{b}_{k}k^{q}$ and
calculate $\hat{b}_{n}=\hat{B}n^{-q}$.
\end{tabular}

\end{center}

\section{Estimation of the Entire Density in the Presence of a Discontinuity
Point}

\subsection{Density Estimation by Truncated Kernels}

We are typically interested in how the shape of the pdf looks like, as well as
whether it has a discontinuity point. \ Imbens and Lemieux (2008) strongly
recommend graphical analyses in empirical studies on RDD, including
inspections of densities of running variables. \ If the test in the previous
section fails to reject the null of continuity of the pdf $f$ at the cutoff
$c$, the entire density may be re-estimated by the gamma kernel, for example.
\ How should we estimate the entire density if the test rejects the null?

The answer to this question is simple. \ It suffices to compute $\hat{f}%
_{-}\left(  x\right)  $ or $\hat{f}_{+}\left(  x\right)  $ as an estimate of
$f\left(  x\right)  $, depending on the position of the design point $x$.
$\ $To put it in another way, $\hat{f}_{-}\left(  x\right)  $ ($\hat{f}%
_{+}\left(  x\right)  $) can be employed whenever $x<c$ ($x>c$), provided that
$c$ is the only point of discontinuity in $f$, as documented in the theorem
below. \ Although only the bias-variance trade-off is provided there,
asymptotic normality of the estimators can be established similarly to Theorem 1.

\paragraph{\emph{Theorem 2.}}

\textit{Suppose that Assumptions 1-3 hold. \ Then, for }$x>c$\textit{,\ as
}$n\rightarrow\infty$, \textit{\ }%
\begin{align*}
Bias\left\{  \hat{f}_{+}\left(  x\right)  \right\}   & \sim\left\{  f^{\left(
1\right)  }\left(  x\right)  +\frac{x}{2}f^{\left(  2\right)  }\left(
x\right)  \right\}  b,\text{\textit{\ and}}\\
Var\left\{  \hat{f}_{+}\left(  x\right)  \right\}   & \sim\frac{1}{nb^{1/2}%
}\frac{f\left(  x\right)  }{2\sqrt{\pi}x^{1/2}}.
\end{align*}
\textit{On the other hand, for }$x<c$\textit{,\ as }$n\rightarrow\infty$,%
\begin{align*}
Bias\left\{  \hat{f}_{-}\left(  x\right)  \right\}   & \sim\left\{  f^{\left(
1\right)  }\left(  x\right)  +\frac{x}{2}f^{\left(  2\right)  }\right\}
b,\text{\textit{\ and}}\\
Var\left\{  \hat{f}_{-}\left(  x\right)  \right\}   & \sim\left\{
\begin{array}
[c]{ll}%
\frac{1}{nb^{1/2}}\frac{f\left(  x\right)  }{2\sqrt{\pi}x^{1/2}} &
\text{\textit{if }}x/b\rightarrow\infty\\
\frac{1}{nb}\frac{\Gamma\left(  2\kappa+1\right)  }{2^{2\kappa+1}\Gamma
^{2}\left(  \kappa+1\right)  }f\left(  x\right)  & \text{\textit{if }%
}x/b\rightarrow\kappa\in\left(  0,\infty\right)
\end{array}
\right.  .
\end{align*}

Theorem 2 indicates no adversity when $f\left(  x\right)  $ for $x\neq c$ is
estimated by $\hat{f}_{\pm}\left(  x\right)  $. $\ $Observe that $\hat{f}%
_{\pm}\left(  x\right)  $ admit the same bias and variance expansions as the
gamma kernel density estimator $\hat{f}\left(  x\right)  $ does. \ A rationale
is that as the design point $x$ moves away from the truncation point $c$, data
points tend to lie on both sides of $x$ and each truncated kernel is likely to
behave like the gamma kernel. \ We can also see that the variance coefficient
decreases as $x$ increases. \ The shrinking variance coefficient as the design
point $x$ moves away from the origin reflects that more data points can be
pooled to smooth in areas with fewer observations. \ This property is
particularly advantageous to estimating the distributions that have a long
tail with sparse data, such as those of the economic variables mentioned in
Section 1.\ 

\subsection{Convergence Properties of $\hat{f}_{-}\left(  x\right)  $ When the
Density Is Unbounded at the Origin}

Clusterings of observations near the boundary are frequently observed in the
distributions with positive supports. \ In the study of RDD, Figure 1 of
Bertrand, Kamenica and Pan (2015) suggests that the distribution of wives'
relative income within households has a clustering of observations near the
origin, as well as a sharp drop at the point of 1/2 (i.e., the point at which
wives' income shares exceed their husbands'). \ Similarly, in Figure 5 of
McCrary (2008), the distribution of proportion of votes for proposed bills in
the US House of Representatives appears to be unbounded at the boundary of
100\%, as well as a sharp discontinuity at the point of 50\%.\footnote{The
arguments in this section are still valid for this case, if we transform the
original data $X$ to $X^{\prime}:=1-X$ and apply them to the transformed data
$X^{\prime}$.}

The following two theorems document weak consistency and the relative
convergence of $\hat{f}_{-}\left(  x\right)  $ when $f\left(  x\right)  $ is
unbounded at $x=0$. \ 

\paragraph{\emph{Theorem 3.}}

\textit{If }$f\left(  x\right)  $\textit{\ is unbounded at }$x=0$\textit{,
Assumptions 1 holds and }$b+\left(  nb^{2}\right)  ^{-1}\rightarrow
0$\textit{\ as }$n\rightarrow\infty$\textit{, then }$\hat{f}_{-}\left(
0\right)  \overset{p}{\rightarrow}\infty$\textit{.}

\paragraph{\emph{Theorem 4.}}

\textit{Suppose that }$f\left(  x\right)  $\textit{\ is unbounded at }%
$x=0$\textit{\ and continuously differentiable in the neighborhood of the
origin. \ In addition, if Assumption 1 holds and }$b+\left\{  nb^{2}f\left(
x\right)  \right\}  ^{-1}\rightarrow0$\textit{\ as }$n\rightarrow\infty
$\textit{\ and }$x\rightarrow0$\textit{, then }%
\[
\left\vert \frac{\hat{f}_{-}\left(  x\right)  -f\left(  x\right)  }{f\left(
x\right)  }\right\vert \overset{p}{\rightarrow}0\mathit{\ }%
\]
\textit{as }$x\rightarrow0$\textit{.}

It has been demonstrated by Bouezmarni and\thinspace Scaillet (2005) and
Hirukawa and Sakudo (2015) that the weak consistency and relative convergence
for densities unbounded at the origin are peculiar to the density estimators
smoothed by the gamma and generalized gamma kernels. \ The theorems ensure
that $\hat{f}_{-}\left(  x\right)  $ is also a proper estimate for unbounded
densities. \ We can deduce from Theorems 2-4 that all in all, appealing
properties of the gamma kernel density estimator are inherited to $\hat
{f}_{\pm}\left(  x\right)  $. \ \ 

\section{Finite-Sample Performance}

It is widely recognized that asymptotic results on kernel-smoothed tests are
not well transmitted to their finite-sample distributions, which reflects that
omitted terms in the first-order asymptotics on the test statistics are highly
sensitive to their smoothing parameter values in finite samples. \ On the
other hand, there is growing literature that reports nice finite-sample
properties of the estimators and test statistics smoothed by asymmetric
kernels. \ Examples include Kristensen (2010) and Gospodinov and Hirukawa
(2012) for estimation and Fernandes and Grammig (2005), Fernandes, Mendes and
Scaillet (2015), and Hirukawa and Sakudo (2016) for testing. \ To see which
perspective dominates, this section investigates finite-sample performance of
the estimator of the jump-size magnitude and the test statistic for
discontinuity of the density via Monte Carlo simulations.

\subsection{Jump-Size Estimation}

First, we focus on the estimator of the jump-size magnitude $J\left(
c\right)  $. \ As true densities, those of the following two asymmetric
distributions are considered: \
\[%
\begin{array}
[c]{ll}%
\text{1. Gamma:} & f\left(  x\right)  =x^{\alpha-1}\exp\left(  -x/\beta
\right)  \mathbf{1}\left(  x\geq0\right)  /\left\{  \beta^{\alpha}%
\Gamma\left(  \alpha\right)  \right\}  ,\,\left(  \alpha,\beta\right)
=\left(  2.75,1\right)  .\\
\text{2. Weibull:} & f\left(  x\right)  =\left(  \alpha/\beta\right)  \left(
x/\beta\right)  ^{\alpha-1}\exp\left\{  -\left(  x/\beta\right)  ^{\alpha
}\right\}  \mathbf{1}\left(  x\geq0\right)  ,\,\left(  \alpha,\beta\right)
=\left(  1.75,3.5\right)  .
\end{array}
\]
Shapes of these densities can be found in Figure 1. \ For each distribution we
choose two suspected discontinuity points $c$, namely, 30\% quantile
(\textquotedblleft30\%\textquotedblright) and median (\textquotedblleft
Med\textquotedblright); see Table 1 for exact values of the points. \ Because
the gamma and Weibull densities have modes at 1.7500 and 2.1567, respectively,
the two points for each density are located on the left- and right-hand sides
of the mode. \ The sample size is\ $n\in\left\{  500,1000,2000\right\}  $, and
$1,000$ replications are drawn for each combination of the sample size $n$ and
the distribution. \ 

The simulation study compares finite-sample performance of our jump-size
estimator $\tilde{J}\left(  c\right)  $ with McCrary's (2008) BLL estimator
$\hat{J}_{M}\left(  c\right)  $. \ The latter employs the triangular kernel
$K\left(  u\right)  =\left(  1-\left\vert u\right\vert \right)  \mathbf{1}%
\left(  \left\vert u\right\vert \leq1\right)  $, and the bandwidth is chosen
by the method described on p.705 of McCrary (2008). \ For the former, the
smoothing parameter $b$ is selected by the power-optimality criterion for two
test statistics $T_{i}\left(  c\right)  :=\sqrt{nb^{1/2}}\tilde{J}\left(
c\right)  /\sqrt{\tilde{V}_{i}\left(  c\right)  }$ for $i=1,2$, where the
definition of $\tilde{V}_{i}\left(  c\right)  $ is given in Section 2.3.
\ Implementation details are as follows: (i) all critical values in $\hat{\pi
}_{M}\left(  b_{k}\right)  $ are set equal to $z_{0.025}=1.96$; (ii) $\left(
p,q\right)  $ are predetermined by $\left(  p,q\right)  =\left(
1/2,4/9\right)  $; (iii) the interval for $b_{k}$\ is $H_{k}=\left[
0.05,0.50\right]  $; and (iv) three different values of the mixing exponent
$\delta$ are considered, namely, $\delta\in\left\{  0.49,0.64,0.81\right\}  $,
so that the exponents on $\hat{f}_{\pm,b}\left(  c\right)  $ and $\hat{f}%
_{\pm,b/\delta}\left(  c\right)  $ to generate $\tilde{f}_{\pm}\left(
c\right)  $ are $\left(  10/3,-7/3\right)  $, $\left(  5,-4\right)  $ and
$\left(  10,-9\right)  $, respectively. \ 

\begin{center}%
\begin{tabular}
[c]{c}\hline\hline
FIGURE\ 1 AND TABLE 1 ABOUT\ HERE\\\hline\hline
\end{tabular}

\end{center}

Table 1 presents as performance measures the bias (\textquotedblleft
Bias\textquotedblright), standard deviation (\textquotedblleft
StdDev\textquotedblright) and root-mean squared error (\textquotedblleft
RMSE\textquotedblright) of each estimator over 1000 Monte Carlo samples.
\ Since the densities are continuous at $c$ actually, the performance measures
are calculated on the basis of $J\left(  c\right)  =0$. \ Moreover, only the
performance measures with the smoothing parameter $b$ selected for
$T_{2}\left(  c\right)  $ are reported, because there is no substantial
difference between values of $b$ chosen for $T_{1}\left(  c\right)  $ and
$T_{2}\left(  c\right)  $.

It can be immediately found that the RMSE shrinks with the sample size, which
indicates consistency of each estimator. \ Although the Bias of $\hat{J}%
_{M}\left(  c\right)  $ is larger than that of $\tilde{J}\left(  c\right)  $,
the StdDev of the former is smaller, and as a consequence it tends to yield a
smaller RMSE. \ We can also see that $\tilde{J}\left(  c\right)  $ has
extremely small biases for all cases, which confirms that the MBC technique
leads to huge bias reduction. \ The bias-variance trade-off in terms of
$\delta$ within $\tilde{J}\left(  c\right)  $ (as Proposition 2 suggests) can
be also observed.

\subsection{Testing for Discontinuity}

Second, size and power properties of the test statistic $T\left(  c\right)  $
are investigated. \ In what follows, $T_{1}\left(  c\right)  $ and
$T_{2}\left(  c\right)  $ are compared with McCrary's (2008) test statistic
based on the difference between logarithms of two density estimates, denoted
as $T_{M}\left(  c\right)  $. \ Implementation details of each test statistic
are the same as described above. \ The Monte Carlo design in this section is
inspired by Otsu, Xu and Matsushita (2013). \ Let $X$ be drawn with
probability $\gamma$ from the truncated gamma or Weibull distribution with
support on $\left[  0,c\right)  $ and with probability $1-\gamma$ from the one
with support on $\left(  c,\infty\right)  $. \ Unless $\gamma=\Pr\left(  X\leq
c\right)  $, the gamma or Weibull pdf is discontinuous at $c$. \ Also denote
the measure of discontinuity as $d:=\Pr\left(  X\leq c\right)  -\gamma$, where
$d\in\left\{  0.00,0.02,0.04,0.06,0.08,0.10\right\}  $ and $d>0$
($\Leftrightarrow J\left(  c\right)  >0$) suggests a jump of the pdf at $c$.
\ For each statistic, the empirical rejection frequencies of the null
$H_{0}:J\left(  c\right)  =0$ for $d=0$ and $d>0$ indicate its size and power
properties, respectively.

\begin{center}%
\begin{tabular}
[c]{c}\hline\hline
TABLES 2-3 ABOUT\ HERE\\\hline\hline
\end{tabular}

\end{center}

Table 2 presents size properties of $T_{1}\left(  c\right)  $ and
$T_{2}\left(  c\right)  $. \ Each test statistic exhibits mild under-rejection
of the null except a few cases, and the rejection frequencies of $T_{2}\left(
c\right)  $ are closer to the nominal ones. \ The rejection frequencies tend
to decrease with $\delta$, and substantial over-rejection of the null is not
observed for $\delta=0.81$. \ Considering that $\delta=0.81$ also yields
nearly unbiased estimates of $J\left(  c\right)  $, we set $\delta$\ equal to
this value for power comparisons. \ 

Table 3 reports power properties of $T_{1}\left(  c\right)  $ and
$T_{2}\left(  c\right)  $, in comparison with $T_{M}\left(  c\right)  $.
\ Panel (A) refers to the results from the gamma distribution. \ It can be
observed that the rejection frequency of each test statistic for a given $d>0$
approaches to one with the sample size $n$, which indicates consistency of
each test. \ Both $T_{1}\left(  c\right)  $ and $T_{2}\left(  c\right)  $
exhibit good power properties without inflating their sizes, and $T_{2}\left(
c\right)  $ appears to be more powerful than $T_{1}\left(  c\right)  $. \ In
contrast, $T_{M}\left(  c\right)  $ exhibits considerable size distortions,
and nonetheless its power properties look inferior to those of $T_{1}\left(
c\right)  $ and $T_{2}\left(  c\right)  $. \ It may be argued that the gamma
distribution is too advantageous to $T_{1}\left(  c\right)  $ and
$T_{2}\left(  c\right)  $ in that both rely on the gamma kernel. \ Hence, the
simulation study based on the Weibull distribution could be fair, and the
results are reported in Panel (B). \ Indeed, the size properties of
$T_{M}\left(  c\right)  $ are dramatically improved. \ However, it is still
outperformed in terms of power properties by $T_{1}\left(  c\right)  $ and
$T_{2}\left(  c\right)  $. \ Again in this case, it appears that $T_{2}\left(
c\right)  $ has better power properties than $T_{1}\left(  c\right)  $. \ A
possible rationale is that because $\tilde{V}_{2}\left(  c\right)  $ tends to
be smaller than $\tilde{V}_{1}\left(  c\right)  $, as suggested in Proposition
2, $T_{2}\left(  c\right)  $ is likely to have a large value (i.e., tends to
reject the null more often) than $T_{1}\left(  c\right)  $ under the
alternative. \ \ \ 

In sum, Monte Carlo results confirm the following two respects. \ First, the
MBC technique achieves huge bias reduction, and the jump-size estimator
$\tilde{J}\left(  c\right)  $ yields nearly unbiased estimates. \ Second, the
test statistics $T_{1}\left(  c\right)  $ and $T_{2}\left(  c\right)  $
exhibit nice power properties without sacrificing their size properties,
whereas the latter appears to be more powerful than the former. \ It is also
worth emphasizing that the superior performance is based simply on first-order
asymptotic results. \ Therefore, assistance of size-adjusting devices such as
bootstrapping appears to be unnecessary, unlike most of the smoothed tests
employing conventional symmetric kernels.\textbf{\ \ }

\section{Empirical Illustration}

This section applies our estimation and testing procedures of discontinuity in
densities to real data. \ We employ the data sets on fourth and fifth graders
of Israeli elementary schools used by Angrist and Lavy (1999). \ The data sets
are made public on the Angrist Data Archive web page, and they are often
utilized in empirical application parts of the closely related literature
(e.g., Otsu, Xu and Matsushita, 2013; Feir,\thinspace\thinspace Lemieux and
Marmer, 2016). \ 

Following Maimonides' rule, Israeli public schools make each class size no
greater than 40. \ As a result of strategic behavior on schools' and/or
parents' sides, the density of school enrollment counts for each grade may be
discontinuous at multiples of 40. \ Then, setting the cutoff $c=40,80,120,160$
for enrollment densities of fourth and fifth graders, we estimate the jump
size and conduct the test for the null of continuity at each cutoff.
\ Specifically, the results from our truncated gamma-kernel approach are
compared with those from McCrary's (2008) BLL method. \ $T_{2}\left(
c\right)  $ with $\delta=0.81$\ is chosen as our test statistic because of its
better finite-sample properties. \ The smoothing parameter for our approach
and the bandwidth for the BLL method are chosen in the same manners as in
Section 4.

\begin{center}%
\begin{tabular}
[c]{c}\hline\hline
FIGURE\ 2 AND TABLE 4 ABOUT\ HERE\\\hline\hline
\end{tabular}

\end{center}

Table 4 presents estimation and testing results on discontinuity in enrollment
densities, where $\hat{f}_{-}^{M}\left(  c\right)  $ and $\hat{f}_{+}%
^{M}\left(  c\right)  $ are BLL estimates of left and right limits of the
density at the cutoff $c$, respectively. \ For convenience, density estimates
with (possible) discontinuity points at $c=40,120$ are plotted in Figure 2.
\ Table 4 shows remarkable differences between estimation results from
McCrary's (2008) and our procedures. \ The former finds upward jump estimates
only at $c=40$ for each grade. \ On the other hand, the latter yields upward
jump estimates at $c=40,80$ and downward jump (or drop) estimates at
$c=120,160$ for each grade. \ In addition, Figure 2 illustrates that the
truncated gamma density estimators tend to capture peaks and troughs more
clearly. \ Testing results also differ. \ While McCrary's (2008) test rejects
the null of continuity at the cutoff only for three cases (i.e., $c=40,120$
for fourth graders and $c=40$ for fifth graders), rejections of the null by
our test include additional two cases (i.e., $c=160 $ for fourth graders and
$c=120$ for fifth graders) as well as the three cases. \ This appears to
reflect better finite-sample power properties of $T_{2}\left(  c\right)  $
reported in Section 4.

\section{Conclusion}

This paper has developed estimation and testing procedures on discontinuity in
densities with positive support. \ Our proposal is built on smoothing by the
gamma kernel. \ To preserve its appealing properties, we split the gamma
kernel into two parts at a given (dis)continuity point and construct two
truncated kernels after re-normalization. \ The jump-size magnitude of the
density at the point can be estimated nonparametrically by two truncated
kernels and the MBC technique by Terrell and Scott (1980). \ The estimator is
easy to implement, and its convergence properties are explored by means of
various approximation techniques on incomplete gamma functions. \ Given the
jump-size estimator, two versions of test statistics for the null of
continuity at a given point are also proposed, and a smoothing parameter
selection method under the power-optimality criterion is tailored to our
testing procedure. \ Furthermore, estimation theory of the entire density in
the presence of a discontinuity point is provided. \ It is demonstrated that
density estimators smoothed by the truncated gamma kernels admit the same bias
and variance approximations as the gamma kernel density estimator does.
\ Monte Carlo simulations indicate that the jump-size estimator is nearly
unbiased when there is no jump in the true density, and that the test
statistics with power-optimal smoothing parameter values plugged in enjoy more
power than McCrary's (2008) BLL-based test does, without sacrificing their
size properties.

We conclude this paper by noting a few research extensions. \ First, the
assumption of a single (known) point of discontinuity may be relaxed. \ It is
worth investigating the cases for more than one (known) point of discontinuity
or those for even unknown (finite) number of discontinuity points. \ For the
latter, locations of discontinuity points are estimated first and then the
corresponding upper and lower limits of the density can be evaluated at each
estimated location. \ Second, while our focus has been exclusively on
univariate densities, the discontinuity analysis may be extended to
multivariate densities. \ 

\appendix\renewcommand{\theequation}{A\arabic{equation}}\setcounter{equation}{0}

\section{Appendix}

\subsection{List of Useful Formulae}

The formulae below are frequently used in the technical proofs.

\begin{description}
\item[\textit{Stirling's formula.}]
\begin{equation}
\Gamma\left(  a+1\right)  =\sqrt{2\pi}a^{a+1/2}\exp\left(  -a\right)  \left\{
1+\frac{1}{12a}+O\left(  a^{-2}\right)  \right\}  \text{ as }a\rightarrow
\infty.\label{UF1}%
\end{equation}

\item[\textit{Recursive formulae on incomplete gamma functions.}]
\begin{align}
\gamma\left(  a+1,z\right)   & =a\gamma\left(  a,z\right)  -z^{a}\exp\left(
-z\right)  \text{ for }a,z>0.\label{UF2}\\
\Gamma\left(  a+1,z\right)   & =a\Gamma\left(  a,z\right)  +z^{a}\exp\left(
-z\right)  \text{ for }a,z>0.\label{UF3}%
\end{align}

\item[\textit{Identity among gamma and incomplete gamma functions. }]
\begin{equation}
\gamma\left(  a,z\right)  +\Gamma\left(  a,z\right)  =\Gamma\left(  a\right)
\text{ for }a,z>0.\label{UF4}%
\end{equation}

\end{description}

\subsection{Proof of Proposition 1}

To save space, we only provide approximations to the bias and variance of
$\hat{f}_{-}\left(  c\right)  $. \ Using (\ref{UF3}), (\ref{UF4}) and
(\ref{P1-1}) gives the results for $\hat{f}_{+}\left(  c\right)  $ in the same
manner. \ The proof utilizes the following asymptotic expansion:%
\begin{equation}
\frac{\gamma\left(  a,a\right)  }{\Gamma\left(  a\right)  }=\frac{1}{2}%
+\frac{1}{\sqrt{2\pi}}\left\{  \frac{1}{3a^{1/2}}+\frac{1}{540a^{3/2}%
}+O\left(  a^{-5/2}\right)  \right\}  \text{ as }a\rightarrow\infty
.\label{P1-1}%
\end{equation}
This can be obtained by either letting $x\downarrow0$ in equation (1) of
Pagurova (1965) or putting $\eta=0$ in equation (1.4) of Temme (1979). \ Then,
putting $z=a$ in (\ref{UF2}) and then substituting (\ref{UF1}) and
(\ref{P1-1}), we have%
\begin{align}
\frac{\gamma\left(  a+1,a\right)  }{\Gamma\left(  a+1\right)  }  &
=\frac{\gamma\left(  a,a\right)  }{\Gamma\left(  a\right)  }-\frac{a^{a}%
\exp\left(  -a\right)  }{\Gamma\left(  a+1\right)  }\nonumber\\
& =\frac{1}{2}+\frac{1}{\sqrt{2\pi}}\left(  -\frac{2}{3}a^{-1/2}+\frac
{23}{270}a^{-3/2}\right)  +O\left(  a^{-5/2}\right)  .\label{P1-2}%
\end{align}

\paragraph{Bias.}

By the change of variable $v:=u/b$,
\[
E\left\{  \hat{f}_{-}\left(  c\right)  \right\}  =\int_{0}^{c}\frac
{u^{c/b}\exp\left(  -u/b\right)  }{b^{c/b+1}\gamma\left(  c/b+1,c/b\right)
}f\left(  u\right)  du=\int_{0}^{a}f\left(  bv\right)  \left\{  \frac
{v^{a}\exp\left(  -v\right)  }{\gamma\left(  a+1,a\right)  }\right\}  dv,
\]
where $a:=c/b$ and the object inside brackets of the right-hand side is a pdf
on the interval $\left[  0,a\right]  $. \ Then, a second-order\textbf{\ }%
Taylor expansion of $f\left(  bv\right)  $ around $bv=c$ (from below) yields
\begin{align}
E\left\{  \hat{f}_{-}\left(  c\right)  \right\}   & =f_{-}\left(  c\right)
+bf_{-}^{\left(  1\right)  }\left(  c\right)  \left\{  \frac{\gamma\left(
a+2,a\right)  }{\gamma\left(  a+1,a\right)  }-a\right\} \nonumber\\
& +\frac{b^{2}}{2}f_{-}^{\left(  2\right)  }\left(  c\right)  \left\{
\frac{\gamma\left(  a+3,a\right)  }{\gamma\left(  a+1,a\right)  }%
-2a\frac{\gamma\left(  a+2,a\right)  }{\gamma\left(  a+1,a\right)  }%
+a^{2}\right\}  +R_{\hat{f}_{-}\left(  c\right)  },\label{P1-3}%
\end{align}
where
\[
R_{\hat{f}_{-}\left(  c\right)  }:=\frac{b^{2}}{2}\int_{0}^{a}\left\{
f_{-}^{\left(  2\right)  }\left(  \xi\right)  -f_{-}^{\left(  2\right)
}\left(  c\right)  \right\}  \left(  v-a\right)  ^{2}\left\{  \frac{v^{a}%
\exp\left(  -v\right)  }{\gamma\left(  a+1,a\right)  }\right\}  dv
\]
is the remainder term with $\xi=\theta\left(  bv\right)  +\left(
1-\theta\right)  c$ for some $\theta\in\left(  0,1\right)  $. \ 

We approximate the leading bias terms first. \ Using (\ref{UF2}) recursively,
we have%
\begin{align*}
\gamma\left(  a+2,a\right)   & =\left(  a+1\right)  \gamma\left(
a+1,a\right)  -a^{a+1}\exp\left(  -a\right)  ,\text{ and}\\
\gamma\left(  a+3,a\right)   & =\left(  a+2\right)  \left(  a+1\right)
\gamma\left(  a+1,a\right)  -2\left(  a+1\right)  a^{a+1}\exp\left(
-a\right)  .
\end{align*}
It follows from (\ref{UF1}) and (\ref{P1-2}) that%
\begin{align*}
\frac{\gamma\left(  a+2,a\right)  }{\gamma\left(  a+1,a\right)  }-a  &
=1-\frac{a^{a+1}\exp\left(  -a\right)  }{\Gamma\left(  a+1\right)  }\left\{
\frac{\gamma\left(  a+1,a\right)  }{\Gamma\left(  a+1\right)  }\right\}
^{-1}\\
& =-\sqrt{\frac{2}{\pi}}a^{1/2}+\left(  1-\frac{4}{3\pi}\right)  +O\left(
a^{-1/2}\right)  ,\text{ and}\\
\frac{\gamma\left(  a+3,a\right)  }{\gamma\left(  a+1,a\right)  }%
-2a\frac{\gamma\left(  a+2,a\right)  }{\gamma\left(  a+1,a\right)  }+a^{2}  &
=a+2-2\frac{a^{a+1}\exp\left(  -a\right)  }{\Gamma\left(  a+1\right)
}\left\{  \frac{\gamma\left(  a+1,a\right)  }{\Gamma\left(  a+1\right)
}\right\}  ^{-1}\\
& =a+O\left(  a^{1/2}\right)  .
\end{align*}
Substituting these into the second and third terms on the right-hand side of
(\ref{P1-3}) and recognizing that $a=c/b$, we obtain%
\begin{align*}
& bf_{-}^{\left(  1\right)  }\left(  c\right)  \left\{  \frac{\gamma\left(
a+2,a\right)  }{\gamma\left(  a+1,a\right)  }-a\right\}  +\frac{b^{2}}{2}%
f_{-}^{\left(  2\right)  }\left(  c\right)  \left\{  \frac{\gamma\left(
a+3,a\right)  }{\gamma\left(  a+1,a\right)  }-2a\frac{\gamma\left(
a+2,a\right)  }{\gamma\left(  a+1,a\right)  }+a^{2}\right\} \\
& =-\sqrt{\frac{2}{\pi}}c^{1/2}f_{-}^{\left(  1\right)  }\left(  c\right)
b^{1/2}+\left\{  \left(  1-\frac{4}{3\pi}\right)  f_{-}^{\left(  1\right)
}\left(  c\right)  +\frac{c}{2}f_{-}^{\left(  2\right)  }\left(  c\right)
\right\}  b+o\left(  b\right)  .
\end{align*}

The remaining task is to demonstrate that $R_{\hat{f}_{-}\left(  c\right)
}=o\left(  b\right)  $. \ It follows from H\"{o}lder-continuity of $f^{\left(
2\right)  }\left(  \cdot\right)  $ and $v\leq c/b=a$ that
\[
\left\vert f^{\left(  2\right)  }\left(  \xi\right)  -f_{-}^{\left(  2\right)
}\left(  c\right)  \right\vert \leq L\left\vert \xi-c\right\vert ^{\varsigma
}=L\theta^{\varsigma}b^{\varsigma}\left(  a-v\right)  ^{\varsigma}.
\]
Using H\"{o}lder's inequality and the fact that $v^{a}\exp\left(  -v\right)
/\gamma\left(  a+1,a\right)  $ is a density on $\left[  0,a\right]  $, we
have
\begin{align*}
\left\vert R_{\hat{f}_{-}\left(  c\right)  }\right\vert  & \leq\frac
{L\theta^{\varsigma}}{2}b^{2+\varsigma}\int_{0}^{a}\left(  a-v\right)
^{2+\varsigma}\left\{  \frac{v^{a}\exp\left(  -v\right)  }{\gamma\left(
a+1,a\right)  }\right\}  dv\\
& \leq\frac{L\theta^{\varsigma}}{2}b^{2+\varsigma}\left[  \int_{0}^{a}\left(
a-v\right)  ^{3}\left\{  \frac{v^{a}\exp\left(  -v\right)  }{\gamma\left(
a+1,a\right)  }\right\}  dv\right]  ^{\left(  2+\varsigma\right)  /3},
\end{align*}
where
\begin{align*}
\int_{0}^{a}\left(  a-v\right)  ^{3}\left\{  \frac{v^{a}\exp\left(  -v\right)
}{\gamma\left(  a+1,a\right)  }\right\}  dv  & =a^{3}-3a^{2}\frac
{\gamma\left(  a+2,a\right)  }{\gamma\left(  a+1,a\right)  }+3a\frac
{\gamma\left(  a+3,a\right)  }{\gamma\left(  a+1,a\right)  }-\frac
{\gamma\left(  a+4,a\right)  }{\gamma\left(  a+1,a\right)  }\\
& =O\left(  a^{3/2}\right)
\end{align*}
by using (\ref{UF1}) and (\ref{P1-2}) repeatedly. \ Finally, substituting
$a=c/b$ yields%
\[
\left\vert R_{\hat{f}_{-}\left(  c\right)  }\right\vert \leq O\left(
b^{2+\varsigma}\right)  O\left\{  b^{-\left(  1+\varsigma/2\right)  }\right\}
=O\left(  b^{1+\varsigma/2}\right)  =o\left(  b\right)  ,
\]
which establishes the bias approximation.

\paragraph{Variance.}

In%
\[
Var\left\{  \hat{f}_{-}\left(  c\right)  \right\}  =\frac{1}{n}E\left\{
K_{G\left(  c,b;c\right)  }^{-}\left(  X_{i}\right)  \right\}  ^{2}+O\left(
n^{-1}\right)  ,
\]
we make an approximation to $E\left\{  K_{G\left(  c,b;c\right)  }^{-}\left(
X_{i}\right)  \right\}  ^{2}$. \ By the change of variable $w:=2u/b$ and
$a=c/b $,
\begin{align*}
E\left\{  K_{G\left(  c,b;c\right)  }^{-}\left(  X_{i}\right)  \right\}  ^{2}
& =\int_{0}^{c}\frac{u^{2c/b}\exp\left(  -2u/b\right)  }{b^{2\left(
c/b+1\right)  }\gamma^{2}\left(  c/b+1,c/b\right)  }f\left(  u\right)  du\\
& =b^{-1}\frac{\gamma\left(  2a+1,2a\right)  }{2^{2a+1}\gamma^{2}\left(
a+1,a\right)  }\int_{0}^{2a}f\left(  \frac{bw}{2}\right)  \left\{
\frac{w^{2a}\exp\left(  -w\right)  }{\gamma\left(  2a+1,2a\right)  }\right\}
dw,
\end{align*}
where the object inside brackets of the right-hand side is again a pdf. \ As
before, the integral part can be approximated by $f_{-}\left(  c\right)
+O\left(  b^{1/2}\right)  $. \ Moreover, it follows from (\ref{P1-2}), the
argument on p.474 of Chen (2000) and $a=c/b$\ that the multiplier part is%
\[
\left\{  \frac{\gamma\left(  2a+1,2a\right)  }{\Gamma\left(  2a+1\right)
}\right\}  \left\{  \frac{\gamma\left(  a+1,a\right)  }{\Gamma\left(
a+1\right)  }\right\}  ^{-2}\left\{  \frac{b^{-1}\Gamma\left(  2a+1\right)
}{2^{2a+1}\Gamma^{2}\left(  a+1\right)  }\right\}  =\frac{b^{-1/2}}{\sqrt{\pi
}c^{1/2}}+o\left(  b^{-1/2}\right)  .
\]
Therefore,%
\[
Var\left\{  \hat{f}_{-}\left(  c\right)  \right\}  =\frac{1}{nb^{1/2}}%
\frac{f_{-}\left(  c\right)  }{\sqrt{\pi}c^{1/2}}+o\left(  n^{-1}%
b^{-1/2}\right)  .\ \blacksquare
\]
\ 

\subsection{Proof of Theorem 1}

The proof requires the following lemma.

\paragraph{\emph{Lemma A1.}}%

\[
E\left\{  K_{G\left(  c,b;c\right)  }^{\pm}\left(  X_{i}\right)  \right\}
^{3}=O\left(  b^{-1}\right)  .
\]

\subsubsection{Proof of Lemma A1}

To save space, we concentrate only on $E\left\{  K_{G\left(  c,b;c\right)
}^{-}\left(  X_{i}\right)  \right\}  ^{3}$. \ By the change of variable
$t:=3u/b $ and $a=c/b$,
\begin{align*}
E\left\{  K_{G\left(  c,b;c\right)  }^{-}\left(  X_{i}\right)  \right\}  ^{3}
& =\int_{0}^{c}\frac{u^{3c/b}\exp\left(  -3u/b\right)  }{b^{3\left(
c/b+1\right)  }\gamma^{3}\left(  c/b+1,c/b\right)  }f\left(  u\right)  du\\
& =b^{-2}\frac{\gamma\left(  3a+1,3a\right)  }{3^{3a+1}\gamma^{3}\left(
a+1,a\right)  }\int_{0}^{3a}f\left(  \frac{bt}{3}\right)  \left\{
\frac{t^{3a}\exp\left(  -t\right)  }{\gamma\left(  3a+1,3a\right)  }\right\}
dt,
\end{align*}
where the integral part is $f_{-}\left(  c\right)  +O\left(  b^{1/2}\right)  $
as before. \ On the other hand, by (\ref{UF1}) and (\ref{P1-2}), the
multiplier part can be approximated by%
\[
\left\{  \frac{\gamma\left(  3a+1,3a\right)  }{\Gamma\left(  3a+1\right)
}\right\}  \left\{  \frac{\gamma\left(  a+1,a\right)  }{\Gamma\left(
a+1\right)  }\right\}  ^{-3}\left\{  \frac{b^{-2}\Gamma\left(  3a+1\right)
}{3^{3a+1}\Gamma^{3}\left(  a+1\right)  }\right\}  =\frac{2}{\sqrt{3}\pi
c}b^{-1}+o\left(  b^{-1}\right)  ,
\]
which establishes the stated result. \ $\blacksquare$

\subsubsection{Proof of Theorem 1}

Let
\begin{align*}
\hat{f}_{\pm,b}\left(  c\right)   & =E\left\{  \hat{f}_{\pm,b}\left(
c\right)  \right\}  +\left[  \hat{f}_{\pm,b}\left(  c\right)  -E\left\{
\hat{f}_{\pm,b}\left(  c\right)  \right\}  \right]  :=I_{b}^{\pm}\left(
c\right)  +Z^{\pm}\text{, and }\\
\hat{f}_{\pm,b/\delta}\left(  c\right)   & =E\left\{  \hat{f}_{\pm,b/\delta
}\left(  c\right)  \right\}  +\left[  \hat{f}_{\pm,b/\delta}\left(  c\right)
-E\left\{  \hat{f}_{\pm,b/\delta}\left(  c\right)  \right\}  \right]
:=I_{b/\delta}^{\pm}\left(  c\right)  +W^{\pm}.
\end{align*}
Then, by a similar argument to the proof for Theorem 1 of Hirukawa and Sakudo
(2014) and Proposition 2,
\begin{align*}
\tilde{J}\left(  c\right)   & =\left\{  I_{b}^{+}\left(  c\right)  \right\}
^{\frac{1}{1-\delta^{1/2}}}\left\{  I_{b/\delta}^{+}\left(  c\right)
\right\}  ^{-\frac{\delta^{1/2}}{1-\delta^{1/2}}}-\left\{  I_{b}^{-}\left(
c\right)  \right\}  ^{\frac{1}{1-\delta^{1/2}}}\left\{  I_{b/\delta}%
^{-}\left(  c\right)  \right\}  ^{-\frac{\delta^{1/2}}{1-\delta^{1/2}}}\\
& +\left(  \frac{1}{1-\delta^{1/2}}\right)  \left\{  \left(  Z^{+}%
-\delta^{1/2}W^{+}\right)  -\left(  Z^{-}-\delta^{1/2}W^{-}\right)  \right\}
+R_{\tilde{J}\left(  c\right)  },
\end{align*}
where it can be shown that the remainder term $R_{\tilde{J}\left(  c\right)
}=o_{p}\left(  n^{-1/2}b^{-1/4}\right)  $. \ Because $E\left(  Z^{\pm}\right)
=E\left(  W^{\pm}\right)  =0$,
\begin{align*}
E\left\{  \tilde{J}\left(  c\right)  \right\}   & \sim\left\{  I_{b}%
^{+}\left(  c\right)  \right\}  ^{\frac{1}{1-\delta^{1/2}}}\left\{
I_{b/\delta}^{+}\left(  c\right)  \right\}  ^{-\frac{\delta^{1/2}}%
{1-\delta^{1/2}}}-\left\{  I_{b}^{-}\left(  c\right)  \right\}  ^{\frac
{1}{1-\delta^{1/2}}}\left\{  I_{b/\delta}^{-}\left(  c\right)  \right\}
^{-\frac{\delta^{1/2}}{1-\delta^{1/2}}}\\
& \sim J\left(  c\right)  +B\left(  c\right)  b,
\end{align*}
where
\begin{align*}
B\left(  c\right)   & =\left(  \frac{1}{\delta^{1/2}}\right)  \left[  \frac
{c}{\pi}\left\{  \frac{\left(  f_{+}^{\left(  1\right)  }\left(  c\right)
\right)  ^{2}}{f_{+}\left(  c\right)  }-\frac{\left(  f_{-}^{\left(  1\right)
}\left(  c\right)  \right)  ^{2}}{f_{-}\left(  c\right)  }\right\}  \right. \\
& \left.  -\left\{  \left(  1-\frac{4}{3\pi}\right)  \left(  f_{+}^{\left(
1\right)  }\left(  c\right)  -f_{-}^{\left(  1\right)  }\left(  c\right)
\right)  +\frac{c}{2}\left(  f_{+}^{\left(  2\right)  }\left(  c\right)
-f_{-}^{\left(  2\right)  }\left(  c\right)  \right)  \right\}  \right]  .
\end{align*}
Therefore,%
\begin{align*}
\sqrt{nb^{1/2}}\left\{  \tilde{J}\left(  c\right)  -J\left(  c\right)
\right\}   & =\sqrt{nb^{1/2}}\left[  \tilde{J}\left(  c\right)  -E\left\{
\tilde{J}\left(  c\right)  \right\}  \right]  +\sqrt{nb^{1/2}}\left[
E\left\{  \tilde{J}\left(  c\right)  \right\}  -J\left(  c\right)  \right] \\
& =\sqrt{nb^{1/2}}\left(  \frac{1}{1-\delta^{1/2}}\right)  \left\{  \left(
Z^{+}-\delta^{1/2}W^{+}\right)  -\left(  Z^{-}-\delta^{1/2}W^{-}\right)
\right\} \\
& +\sqrt{nb^{1/2}}\left\{  B\left(  c\right)  b+o\left(  b\right)  \right\}
+o_{p}\left(  1\right)  ,
\end{align*}
where the second term on the right hand side becomes asymptotically negligible
if $nb^{5/2}\rightarrow0$. \ 

The remaining task is to establish the asymptotic normality of the first term.
\ Due to the disjunction of two truncated kernels $K_{G\left(  c,b;c\right)
}^{\pm}\left(  \cdot\right)  $, the asymptotic variance of the term, denoted
as $V\left(  c\right)  $, is just the sum of asymptotic variances of
$\tilde{f}_{\pm}\left(  c\right)  $ given in Proposition 2. \ Hence, we need
only to establish Liapunov's condition. \ Denoting%
\begin{align*}
Z^{\pm}  & =\sum_{i=1}^{n}\left(  \frac{1}{n}\right)  \left[  K_{G\left(
c,b;c\right)  }^{\pm}\left(  X_{i}\right)  -E\left\{  K_{G\left(
c,b;c\right)  }^{\pm}\left(  X_{i}\right)  \right\}  \right]  :=\sum_{i=1}%
^{n}\left(  \frac{1}{n}\right)  Z_{i}^{\pm},\text{ and}\\
W^{\pm}  & =\sum_{i=1}^{n}\left(  \frac{1}{n}\right)  \left[  K_{G\left(
c,b/\delta;c\right)  }^{\pm}\left(  X_{i}\right)  -E\left\{  K_{G\left(
c,b/\delta;c\right)  }^{\pm}\left(  X_{i}\right)  \right\}  \right]
:=\sum_{i=1}^{n}\left(  \frac{1}{n}\right)  W_{i}^{\pm},
\end{align*}
we can rewrite the term as%
\begin{align*}
& \sqrt{nb^{1/2}}\left(  \frac{1}{1-\delta^{1/2}}\right)  \left\{  \left(
Z^{+}-\delta^{1/2}W^{+}\right)  -\left(  Z^{-}-\delta^{1/2}W^{-}\right)
\right\} \\
& =\sum_{i=1}^{n}\sqrt{\frac{b^{1/2}}{n}}\left(  \frac{1}{1-\delta^{1/2}%
}\right)  \left\{  \left(  Z_{i}^{+}-\delta^{1/2}W_{i}^{+}\right)  -\left(
Z_{i}^{-}-\delta^{1/2}W_{i}^{-}\right)  \right\}  :=\sum_{i=1}^{n}Y_{i}.
\end{align*}
It follows from $0<\delta<1$ that%
\[
E\left\vert Y_{i}\right\vert ^{3}\leq\frac{b^{3/4}}{n^{3/2}}\left(  \frac
{1}{1-\delta^{1/2}}\right)  ^{3}E\left(  \left\vert Z_{i}^{+}\right\vert
+\left\vert W_{i}^{+}\right\vert +\left\vert Z_{i}^{-}\right\vert +\left\vert
W_{i}^{-}\right\vert \right)  ^{3}.
\]
Because the expected value part is $O\left(  b^{-1}\right)  $ by Lemma A1,
$E\left\vert Y_{i}\right\vert ^{3}=O\left(  n^{-3/2}b^{-1/4}\right)  $. \ It
is also straightforward to see that $Var\left(  Y_{i}\right)  =O\left(
n^{-1}\right)  $. \ Therefore,
\[
\frac{\sum_{i=1}^{n}E\left\vert Y_{i}\right\vert ^{3}}{\left\{  \sum_{i=1}%
^{n}Var\left(  Y_{i}\right)  \right\}  ^{3/2}}=O\left(  n^{-1/2}%
b^{-1/4}\right)  \rightarrow0,
\]
or Liapunov's condition holds. \ This completes the proof. \ $\blacksquare$

\subsection{Proof of Proposition 3}

The proof closely follows the one for Proposition 1 of Hirukawa and Sakudo
(2016). \ It follows from Theorem 1 that $E\left\{  \tilde{J}\left(  c\right)
\right\}  =J\left(  c\right)  +O\left(  b\right)  $, $Var\left\{  \tilde
{J}\left(  c\right)  \right\}  =O\left(  n^{-1}b^{-1/2}\right)  $ and
$\tilde{V}\left(  c\right)  \overset{p}{\rightarrow}V\left(  c\right)  $,
regardless of whether $H_{0}$ or $H_{1}$ may be true. \ Therefore, $\tilde
{J}\left(  c\right)  =J\left(  c\right)  +O\left(  b\right)  +O_{p}\left(
n^{-1/2}b^{-1/4}\right)  \overset{p}{\rightarrow}J\left(  c\right)  \neq0$
under $H_{1}$, and thus $\left\vert T\left(  c\right)  \right\vert $ is a
divergent stochastic sequence with an expansion rate of $n^{1/2}b^{1/4}$.
\ The result immediately follows. $\ \blacksquare$

\subsection{Proof of Theorem 2}

To demonstrate this theorem, we must rely on different asymptotic expansions,
depending on the positions of the design point $x$ and the truncation point
$c$. \ For notational convenience, put $\left(  a,z\right)  =\left(
x/b,c/b\right)  $. \ The proof requires the following lemma.

\paragraph{\emph{Lemma A2.}}

\textit{For }$a>0$\textit{\ and }$z>\max\left\{  1,a\right\}  $\textit{, }%
\[
\Gamma\left(  a+1,z\right)  \leq\left\{
\begin{array}
[c]{ll}%
z^{a}\exp\left(  -z\right)  +\exp\left(  -z\right)  & \text{\textit{for }%
}0<a\leq1\\
\left(  a+1\right)  z^{a}\exp\left(  -z\right)  +\Gamma\left(  a+1\right)
\exp\left(  -z\right)  & \text{\textit{for }}a>1
\end{array}
\right.  .
\]

\subsubsection{Proof of Lemma A2}

For $0<a\leq1$, it follows from an elementary inequality on the upper
incomplete gamma function (e.g., equation (1.05) on p.67 of Olver, 1974) and
$z>1$ that
\begin{equation}
\Gamma\left(  a,z\right)  \leq z^{a-1}\exp\left(  -z\right)  \leq\exp\left(
-z\right)  .\label{LA2-1}%
\end{equation}
Then, by (\ref{UF3}),%
\[
\Gamma\left(  a+1,z\right)  =z^{a}\exp\left(  -z\right)  +a\Gamma\left(
a,z\right)  \leq z^{a}\exp\left(  -z\right)  +1\cdot\exp\left(  -z\right)  .
\]
Next, for $a>1$ and $a\in%
\mathbb{N}
$, using (\ref{UF3}) recursively yields%
\begin{align*}
\Gamma\left(  a+1,z\right)   & =z^{a}\exp\left(  -z\right)  \left\{
1+\frac{a}{z}+\frac{a\left(  a-1\right)  }{z^{2}}+\cdots+\frac{a\left(
a-1\right)  \cdots2}{z^{a-1}}\right\} \\
& +a\left(  a-1\right)  \cdots2\cdot1\cdot\Gamma\left(  1,z\right)  ,
\end{align*}
where the sum inside the brackets is bounded by $a\left(  \leq a+1\right)  $.
\ Then, by (\ref{LA2-1}),%
\[
\Gamma\left(  a+1,z\right)  \leq\left(  a+1\right)  z^{a}\exp\left(
-z\right)  +\Gamma\left(  a+1\right)  \exp\left(  -z\right)  .
\]
Finally, for $a>1$ and $a\notin%
\mathbb{N}
$, we have
\begin{align*}
\Gamma\left(  a+1,z\right)   & =z^{a}\exp\left(  -z\right)  \left\{
1+\frac{a}{z}+\frac{a\left(  a-1\right)  }{z^{2}}+\cdots+\frac{a\left(
a-1\right)  \cdots\left(  a-\left\lfloor a\right\rfloor +1\right)
}{z^{\left\lfloor a\right\rfloor }}\right\} \\
& +a\left(  a-1\right)  \cdots\left(  a-\left\lfloor a\right\rfloor \right)
\Gamma\left(  a-\left\lfloor a\right\rfloor ,z\right)  .
\end{align*}
where the sum inside the brackets is bounded by $\left\lfloor a\right\rfloor
+1\left(  \leq a+1\right)  $. \ Because $0<a-\left\lfloor a\right\rfloor <1$,
$\Gamma\left(  a-\left\lfloor a\right\rfloor \right)  >1$ and thus%
\[
a\left(  a-1\right)  \cdots\left(  a-\left\lfloor a\right\rfloor \right)
=\frac{\Gamma\left(  a+1\right)  }{\Gamma\left(  a-\left\lfloor a\right\rfloor
\right)  }\leq\Gamma\left(  a+1\right)  .
\]
Therefore, again by (\ref{LA2-1}),
\[
\Gamma\left(  a+1,z\right)  \leq\left(  a+1\right)  z^{a}\exp\left(
-z\right)  +\Gamma\left(  a+1\right)  \exp\left(  -z\right)  .\text{
\ }\blacksquare
\]

\subsubsection{Proof of Theorem 2}

\subsubsection*{(i) \textbf{On }$\hat{f}_{-}\left(  x\right)  $:}

We consider different approximations to incomplete gamma functions depending
on the position of $x$. \ When $x/b\rightarrow\infty$, $z>a$ and
$a,z\rightarrow\infty$ hold. \ Hence, the case for $a>1$ of Lemma A2 applies,
and thus%
\[
\frac{\Gamma\left(  a+1,z\right)  }{\Gamma\left(  a+1\right)  }\leq\left(
a+1\right)  \left\{  \frac{z^{a}\exp\left(  -z\right)  }{\Gamma\left(
a+1\right)  }\right\}  +\exp\left(  -z\right)  .
\]
It follows from (\ref{UF1}) and $\rho:=a/z\in\left(  0,1\right)  $ that%
\begin{align}
\frac{z^{a}\exp\left(  -z\right)  }{\Gamma\left(  a+1\right)  }  & =\left\{
\frac{1+O\left(  a^{-1}\right)  }{\sqrt{2\pi}}\right\}  a^{-1/2}\exp\left\{
a\ln\left(  \frac{e}{\rho e^{1/\rho}}\right)  \right\} \nonumber\\
& =O\left[  a^{-1/2}\exp\left\{  a\ln\left(  \frac{e}{\rho e^{1/\rho}}\right)
\right\}  \right]  ,\label{Thm2-1}%
\end{align}
where $e/\left(  \rho e^{1/\rho}\right)  \in\left(  0,1\right)  $ holds.
\ Then,%
\[
\frac{\Gamma\left(  a+1,z\right)  }{\Gamma\left(  a+1\right)  }=O\left[
a^{1/2}\exp\left\{  a\ln\left(  \frac{e}{\rho e^{1/\rho}}\right)  \right\}
\right]  .
\]
On the other hand, when $x/b\rightarrow\kappa\in\left(  0,\infty\right)  $,
putting $a\rightarrow\kappa$ and $z\rightarrow\infty$ in Lemma A2 yields%
\[
\frac{\Gamma\left(  a+1,z\right)  }{\Gamma\left(  a+1\right)  }=O\left\{
z^{\kappa}\exp\left(  -z\right)  \right\}  .
\]
By (\ref{UF4}), we finally have
\begin{equation}
\frac{\gamma\left(  a+1,z\right)  }{\Gamma\left(  a+1\right)  }=1+\left\{
\begin{array}
[c]{ll}%
O\left[  a^{1/2}\exp\left\{  a\ln\left(  e/\left(  \rho e^{1/\rho}\right)
\right)  \right\}  \right]  & \text{if }x/b\rightarrow\infty\\
O\left\{  z^{\kappa}\exp\left(  -z\right)  \right\}  & \text{if }%
x/b\rightarrow\kappa
\end{array}
\right.  .\label{Thm2-2}%
\end{equation}

\paragraph{Bias.}

By (\ref{Thm2-1}), (\ref{Thm2-2}), and $\left(  a,z\right)  =\left(
x/b,c/b\right)  $,
\begin{align*}
& \frac{\gamma\left(  a+2,z\right)  }{\gamma\left(  a+1,z\right)  }-a\\
& =1-\frac{z^{a+1}\exp\left(  -z\right)  }{\Gamma\left(  a+1\right)  }\left\{
\frac{\gamma\left(  a+1,z\right)  }{\Gamma\left(  a+1\right)  }\right\}
^{-1}\\
& =1+\left\{
\begin{array}
[c]{ll}%
O\left[  a^{1/2}\exp\left\{  a\ln\left(  e/\left(  \rho e^{1/\rho}\right)
\right)  \right\}  \right]  & \text{if }x/b\rightarrow\infty\\
O\left\{  z^{\kappa}\exp\left(  -z\right)  \right\}  & \text{if }%
x/b\rightarrow\kappa
\end{array}
\right. \\
& =1+\left\{
\begin{array}
[c]{ll}%
O\left[  b^{-1/2}\exp\left\{  \left(  x/b\right)  \ln\left(  e/\left(  \rho
e^{1/\rho}\right)  \right)  \right\}  \right]  & \text{if }x/b\rightarrow
\infty\\
O\left\{  b^{-\kappa}\exp\left(  -c/b\right)  \right\}  & \text{if
}x/b\rightarrow\kappa
\end{array}
\right.  ,\text{ and}\\
& \frac{\gamma\left(  a+3,z\right)  }{\gamma\left(  a+1,z\right)  }%
-2a\frac{\gamma\left(  a+2,z\right)  }{\gamma\left(  a+1,z\right)  }+a^{2}\\
& =a+2-\left(  z-a+2\right)  \frac{z^{a+1}\exp\left(  -z\right)  }%
{\Gamma\left(  a+1\right)  }\left\{  \frac{\gamma\left(  a+1,z\right)
}{\Gamma\left(  a+1\right)  }\right\}  ^{-1}\\
& =a+2+\left\{
\begin{array}
[c]{ll}%
O\left[  a^{3/2}\exp\left\{  a\ln\left(  e/\left(  \rho e^{1/\rho}\right)
\right)  \right\}  \right]  & \text{if }x/b\rightarrow\infty\\
O\left\{  z^{\kappa+1}\exp\left(  -z\right)  \right\}  & \text{if
}x/b\rightarrow\kappa
\end{array}
\right. \\
& =\frac{x}{b}+2+\left\{
\begin{array}
[c]{ll}%
O\left[  b^{-3/2}\exp\left\{  \left(  x/b\right)  \ln\left(  e/\left(  \rho
e^{1/\rho}\right)  \right)  \right\}  \right]  & \text{if }x/b\rightarrow
\infty\\
O\left\{  b^{-\kappa-1}\exp\left(  -c/b\right)  \right\}  & \text{if
}x/b\rightarrow\kappa
\end{array}
\right.  .
\end{align*}
Then, by the argument in the proof of Proposition 1, in either case,
\[
E\left\{  \hat{f}_{-}\left(  x\right)  \right\}  =f\left(  x\right)  +\left\{
f^{\left(  1\right)  }\left(  x\right)  +\frac{x}{2}f^{\left(  2\right)
}\left(  x\right)  \right\}  b+o\left(  b\right)  .
\]

\paragraph{Variance.}

In
\[
E\left\{  K_{G\left(  x,b;c\right)  }^{-}\left(  X_{i}\right)  \right\}
^{2}=b^{-1}\frac{\gamma\left(  2a+1,2z\right)  }{2^{2a+1}\gamma^{2}\left(
a+1,z\right)  }\int_{0}^{2z}f\left(  \frac{bw}{2}\right)  \left\{
\frac{w^{2a}\exp\left(  -w\right)  }{\gamma\left(  2a+1,2z\right)  }\right\}
dw,
\]
the integral part is $f\left(  x\right)  +O\left(  b\right)  $\ in either
case. \ It also follows from (\ref{Thm2-2}) and the argument on p.474 of Chen
(2000) that the multiplier part is%
\begin{align*}
& \left\{  \frac{\gamma\left(  2a+1,2z\right)  }{\Gamma\left(  2a+1\right)
}\right\}  \left\{  \frac{\gamma\left(  a+1,z\right)  }{\Gamma\left(
a+1\right)  }\right\}  ^{-2}\left\{  \frac{b^{-1}\Gamma\left(  2a+1\right)
}{2^{2a+1}\Gamma^{2}\left(  a+1\right)  }\right\} \\
& =\left\{
\begin{array}
[c]{ll}%
\frac{b^{-1/2}}{2\sqrt{\pi}x^{1/2}}+o\left(  b^{-1/2}\right)  & \text{if
}x/b\rightarrow\infty\\
\frac{b^{-1}\Gamma\left(  2\kappa+1\right)  }{2^{2\kappa+1}\Gamma^{2}\left(
\kappa+1\right)  }+o\left(  b^{-1}\right)  & \text{if }x/b\rightarrow\kappa
\end{array}
\right.  .
\end{align*}
Therefore,%
\[
Var\left\{  \hat{f}_{-}\left(  x\right)  \right\}  =\left\{
\begin{array}
[c]{ll}%
\frac{1}{nb^{1/2}}\frac{f\left(  x\right)  }{2\sqrt{\pi}x^{1/2}}+o\left(
n^{-1}b^{-1/2}\right)  & \text{if }x/b\rightarrow\infty\\
\frac{1}{nb}\frac{\Gamma\left(  2\kappa+1\right)  }{2^{2\kappa+1}\Gamma
^{2}\left(  \kappa+1\right)  }f\left(  x\right)  +o\left(  n^{-1}b^{-1}\right)
& \text{if }x/b\rightarrow\kappa
\end{array}
\right.  .\ \blacksquare
\]

\subsubsection*{(ii) \textbf{On }$\hat{f}_{+}\left(  x\right)  $:}

We may focus only on the case for interior $x$. \ However, it seems difficult
to derive a sharp bound on $\gamma\left(  a+1,z\right)  $ or $\Gamma\left(
a+1,z\right)  $ for the case of $a>z$ and $a,z\rightarrow\infty$ based
directly on (\ref{UF2}) or (\ref{UF3}). \ Instead, we turn to the series
expansion described in Section 3 of Ferreira,\thinspace\thinspace L\'{o}pez
and P\'{e}rez-Sinus\'{\i}a (2005), which is valid for the case of $a>z$,
$a,z\rightarrow\infty$ and $a-z=O\left(  a\right)  $. \ The expansion is
\[
\gamma\left(  a+1,z\right)  =z^{a+1}\exp\left(  -z\right)  \sum_{k=0}^{\infty
}c_{k}\left(  a\right)  \Phi_{k}\left(  z-a\right)  ,
\]
where the definitions of $\left\{  c_{k}\left(  a\right)  \right\}  $ and
$\left\{  \Phi_{k}\left(  z-a\right)  \right\}  $ can be found therein.
\ Because the sum is shown to be convergent, the order of magnitude in
$\gamma\left(  a+1,z\right)  /\Gamma\left(  a+1\right)  $ is determined by the
one in $z^{a+1}\exp\left(  -z\right)  /\Gamma\left(  a+1\right)  $. \ It
follows from (\ref{UF1}) and $\rho^{\prime}:=z/a\in\left(  0,1\right)  $ that%
\begin{align*}
\frac{z^{a+1}\exp\left(  -z\right)  }{\Gamma\left(  a+1\right)  }  & =\left[
\frac{\rho^{\prime}\left\{  1+O\left(  a^{-1}\right)  \right\}  }{\sqrt{2\pi}%
}\right]  a^{1/2}\exp\left\{  a\ln\left(  \frac{\rho^{\prime}e}{e^{\rho
^{\prime}}}\right)  \right\} \\
& =O\left[  a^{1/2}\exp\left\{  a\ln\left(  \frac{\rho^{\prime}e}%
{e^{\rho^{\prime}}}\right)  \right\}  \right]  ,
\end{align*}
where $\rho^{\prime}e/e^{\rho^{\prime}}\in\left(  0,1\right)  $ is again the
case. \ Then, by (\ref{UF4}),
\[
\frac{\Gamma\left(  a+1,z\right)  }{\Gamma\left(  a+1\right)  }=1+O\left[
a^{1/2}\exp\left\{  a\ln\left(  \frac{\rho^{\prime}e}{e^{\rho^{\prime}}%
}\right)  \right\}  \right]  .
\]
The bias and variance of $\hat{f}_{+}\left(  x\right)  $\ can be approximated
as above. \ $\blacksquare$

\subsection{Proof of Theorem 3}

Both this proof and the proof of Theorem 4 require three lemmata below. \ 

\paragraph{\emph{Lemma A3.}}

\textit{For }$\alpha>0$\textit{\ and a sufficiently small }$b>0$\textit{, pick
some design point }$x\in\left[  0,\alpha b\right]  $\textit{. \ Then, for
}$\eta\in\left(  0,c\right)  $\textit{, }%
\[
\int_{0}^{\eta}K_{G\left(  x,b;c\right)  }^{-}\left(  u\right)  du=\int
_{0}^{\eta}\frac{u^{x/b}\exp\left(  -u/b\right)  }{b^{x/b+1}\gamma\left(
x/b+1,c/b\right)  }du\rightarrow1
\]
\textit{as }$b\rightarrow0$\textit{.\ \ }

\paragraph{\emph{Lemma A4.}}

\textit{For the design point }$x$\textit{\ defined in Lemma A3, let }%
\[
\left\{  K_{i}\right\}  _{i=1}^{n}:=\left\{  bK_{G\left(  x,b;c\right)  }%
^{-}\left(  X_{i}\right)  \right\}  _{i=1}^{n}.\mathit{\ \ }%
\]
\textit{Then, }%
\[
0\leq K_{i}\leq C:=\max\left\{  1,\alpha^{\alpha}\right\}  \left\{
\frac{\Gamma\left(  \alpha+1\right)  }{\gamma\left(  \alpha+1,\alpha\right)
}\right\}  \left\{  \frac{1}{\Gamma\left(  a^{\ast}\right)  }\right\}  ,
\]
\textit{where }$\Gamma\left(  a^{\ast}\right)  :=\min_{a>0}\Gamma\left(
a\right)  \approx0.8856$\textit{\ for }$a^{\ast}\approx1.4616$\textit{.\ \ }

\paragraph{\emph{Lemma A5.(Hoeffding, 1963, Theorem 2)}}

\textit{Let }$\left\{  X_{i}\right\}  _{i=1}^{n}$\textit{\ be independent and
}$a_{i}\leq X_{i}\leq b_{i}$\textit{\ for }$i=1,2,\ldots,n$\textit{. \ Also
write }$\bar{X}:=\left(  1/n\right)  \sum_{i=1}^{n}X_{i}$\textit{\ and }%
$\mu:=E\left(  \bar{X}\right)  $\textit{. \ Then, for }$\epsilon>0$\textit{,}%
\[
\Pr\left(  \left\vert \bar{X}-\mu\right\vert \geq\epsilon\right)  \leq
2\exp\left\{  -\frac{2n^{2}\epsilon^{2}}{\sum_{i=1}^{n}\left(  b_{i}%
-a_{i}\right)  ^{2}}\right\}  .
\]

\subsubsection{Proof of Lemma A3}

By the change of variable $v:=u/b$, the integral can be rewritten as
\[
\int_{0}^{\eta/b}\frac{v^{x/b}\exp\left(  -v\right)  }{\gamma\left(
x/b+1,c/b\right)  }dv=\frac{\gamma\left(  x/b+1,\eta/b\right)  }{\gamma\left(
x/b+1,c/b\right)  }.
\]
Because $\eta/b\uparrow\infty$ and $0\leq x/b\leq\alpha$, (\ref{Thm2-2})
establishes that \
\[
\frac{\gamma\left(  x/b+1,\eta/b\right)  }{\gamma\left(  x/b+1,c/b\right)
}=\frac{\Gamma\left(  x/b+1\right)  +O\left\{  b^{-\alpha}\exp\left(
-\eta/b\right)  \right\}  }{\Gamma\left(  x/b+1\right)  +O\left\{  b^{-\alpha
}\exp\left(  -c/b\right)  \right\}  }\rightarrow1.\ \blacksquare
\]

\subsubsection{Proof of Lemma A4}

By construction, $K_{i}\geq0$ holds. \ In addition, since the gamma kernel has
its mode at the design point $x$ (Chen, 2000, p.473), $K_{i}$ is bounded by%
\begin{equation}
bK_{G\left(  x,b;c\right)  }^{-}\left(  x\right)  =\left(  \frac{x}{b}\right)
^{x/b}\exp\left(  -\frac{x}{b}\right)  \left\{  \frac{\Gamma\left(
x/b+1\right)  }{\gamma\left(  x/b+1,c/b\right)  }\right\}  \left\{  \frac
{1}{\Gamma\left(  x/b+1\right)  }\right\}  .\label{LA4-1}%
\end{equation}
For $0\leq x/b\leq\alpha$, $\left(  x/b\right)  ^{x/b}\leq\max\left\{
1,\alpha^{\alpha}\right\}  $ and $\exp\left(  -x/b\right)  \leq1$. \ Moreover,
$\gamma\left(  a,z\right)  /\Gamma\left(  a\right)  $ for $a,z>0$ is
monotonously increasing in $z$ and decreasing in $a$; see, for example,
Tricomi (1950, p.276) for details. \ Because $c$ is an interior point, $\alpha
b\leq c$ or $\alpha\leq c/b$ holds. \ Hence,
\[
\frac{\Gamma\left(  x/b+1\right)  }{\gamma\left(  x/b+1,c/b\right)  }\leq
\frac{\Gamma\left(  \alpha+1\right)  }{\gamma\left(  \alpha+1,\alpha\right)
}.\
\]
Finally, it is known that $\Gamma\left(  a^{\ast}\right)  :=\min_{a>0}%
\Gamma\left(  a\right)  \approx0.8856$ for\textit{\ }$a^{\ast}\approx
1.4616$\textit{.} \ Therefore, the right-hand side of (\ref{LA4-1}) has the
upper bound
\[
\max\left\{  1,\alpha^{\alpha}\right\}  \cdot1\cdot\left\{  \frac
{\Gamma\left(  \alpha+1\right)  }{\gamma\left(  \alpha+1,\alpha\right)
}\right\}  \left\{  \frac{1}{\Gamma\left(  a^{\ast}\right)  }\right\}
:=C.\ \blacksquare
\]

\subsubsection{Proof of Theorem 3}

This proof largely follows the one for Theorem 5 of Hirukawa and Sakudo
(2015). \ Without loss of generality, for $\alpha>0$ and a sufficiently small
$b>0$, pick some design point $x\in\left[  0,\alpha b\right]  $. \ Then,\ the
proof completes if the following statements hold:%
\begin{align}
\hat{f}_{-}\left(  x\right)   & =E\left\{  \hat{f}_{-}\left(  x\right)
\right\}  +o_{p}\left(  1\right)  .\label{Thm3-1}\\
E\left\{  \hat{f}_{-}\left(  x\right)  \right\}   & =E\left\{  \hat{f}%
_{-}\left(  0\right)  \right\}  +o\left(  1\right)  .\label{Thm3-2}\\
E\left\{  \hat{f}_{-}\left(  0\right)  \right\}   & \rightarrow\infty
.\label{Thm3-3}%
\end{align}

Below we demonstrate (\ref{Thm3-1})-(\ref{Thm3-3}) one by one. \ First,
(\ref{Thm3-2}) immediately follows from the continuity of $K_{G\left(
x,b;c\right)  }^{-}\left(  u\right)  $ in $x$. \ Second, when $f\left(
x\right)  \rightarrow\infty$ as $x\rightarrow0$, it holds that for any $A>0$,
there is some $\eta\in\left(  0,c\right)  $ such that $f\left(  x\right)  >A$
for all $x<\eta$. \ For the given $\eta$, Lemma A3 implies that
\[
E\left\{  \hat{f}_{-}\left(  0\right)  \right\}  >\int_{0}^{\eta}K_{G\left(
0,b;c\right)  }^{-}\left(  u\right)  f\left(  u\right)  du>A\int_{0}^{\eta
}K_{G\left(  0,b;c\right)  }^{-}\left(  u\right)  du\rightarrow A,
\]
which establishes (\ref{Thm3-3}). \ Third, for $\left\{  K_{i}\right\}
_{i=1}^{n}$ defined in Lemma A4, denote their sample average as $\bar
{K}:=\left(  1/n\right)  \sum_{i=1}^{n}K_{i}$. \ Then, it follows from Lemmata
A4 and A5 that for $\epsilon>0$,%
\begin{align*}
\Pr\left(  \left\vert \hat{f}_{-}\left(  x\right)  -E\left\{  \hat{f}%
_{-}\left(  x\right)  \right\}  \right\vert \geq\epsilon\right)   &
=\Pr\left(  \left\vert \bar{K}-E\left(  K_{i}\right)  \right\vert \geq
b\epsilon\right) \\
& \leq2\exp\left\{  -2\left(  \frac{\epsilon}{C}\right)  ^{2}nb^{2}\right\}
\rightarrow0.
\end{align*}
Therefore, (\ref{Thm3-1})\ is also demonstrated, and thus the proof is
completed. \ $\blacksquare$

\subsection{Proof of Theorem 4}

This proof largely follows the one for Theorem 5.3 of Bouezmarni and Scaillet
(2005). \ As in the proof of Theorem 3, pick some $x\in\left[  0,\alpha
b\right]  $. \ Then, the proof is boiled down to establishing the following
statements:
\begin{align}
& \left\vert \frac{E\left\{  \hat{f}_{-}\left(  x\right)  \right\}  -f\left(
x\right)  }{f\left(  x\right)  }\right\vert \rightarrow0,\text{ and}%
\label{Thm4-1}\\
& \left\vert \frac{\hat{f}_{-}\left(  x\right)  -E\left\{  \hat{f}_{-}\left(
x\right)  \right\}  }{f\left(  x\right)  }\right\vert \overset{p}{\rightarrow
}0,\label{Thm4-2}%
\end{align}
as $n\rightarrow\infty$ and $b,x\rightarrow0$.

We demonstrate (\ref{Thm4-1}) first. \ An inspection of the proof for Theorem
5.3 of Bouezmarni and Scaillet (2005) reveals that (\ref{Thm4-1}) is shown if
their conditions A.2, A.3 and A.5 are fulfilled. \ Now we check the validity
of three conditions. \ First, because $\int_{0}^{\infty}f\left(  x\right)
dx=1$ and $f\left(  x\right)  \rightarrow\infty$ as $x\rightarrow0$, there are
constants $0<\underline{C}<\overline{C}<$ $\infty$ such that $\underline
{C}x^{-d}\leq f\left(  x\right)  \leq\overline{C}x^{-d}$ for some $d\in\left(
0,1\right)  $ as $x\rightarrow0$. \ Accordingly, $f^{\left(  1\right)
}\left(  x\right)  =O\left(  x^{-d-1}\right)  $ for a small value of $x$.
$\ $These imply that $x\left\vert f^{\left(  1\right)  }\left(  x\right)
\right\vert /f\left(  x\right)  \leq O\left(  1\right)  $, and thus A.2
follows. \ Second, A.3 has been already established as Lemma A1. \ Third, let
the random variable $U$ be drawn from the distribution with the pdf
$K_{G\left(  x,b;c\right)  }^{-}\left(  u\right)  $. \ Then, by $0\leq
x/b\leq\alpha$ and the expansion techniques used in the proof of Theorem 2,
$Var\left(  U\right)  \leq O\left(  b\right)  \rightarrow0$, and thus A.5 also holds.\ \ \ 

Furthermore, it follows from Lemmata A4 and A5 that for $\bar{K}$ defined in
the proof of Theorem 3 and for $\epsilon>0$,%
\begin{align*}
\Pr\left(  \left\vert \frac{\hat{f}_{-}\left(  x\right)  -E\left\{  \hat
{f}_{-}\left(  x\right)  \right\}  }{f\left(  x\right)  }\right\vert
\geq\epsilon\right)   & =\Pr\left(  \left\vert \bar{K}-E\left(  K_{i}\right)
\right\vert \geq bf\left(  x\right)  \epsilon\right) \\
& \leq2\exp\left\{  -2\left(  \frac{\epsilon}{C}\right)  ^{2}nb^{2}%
f^{2}\left(  x\right)  \right\}  \rightarrow0.
\end{align*}
Therefore, (\ref{Thm4-2})\ is also demonstrated, and thus the proof is
completed. \ $\blacksquare$

\setlength{\baselineskip}{12pt}

\begin{center}
\bigskip

\pagebreak

\textbf{Figure 1:} Shapes of True Densities for Monte Carlo
Simulations\bigskip%

{\includegraphics[
height=2.2748in,
width=4.6442in
]%
{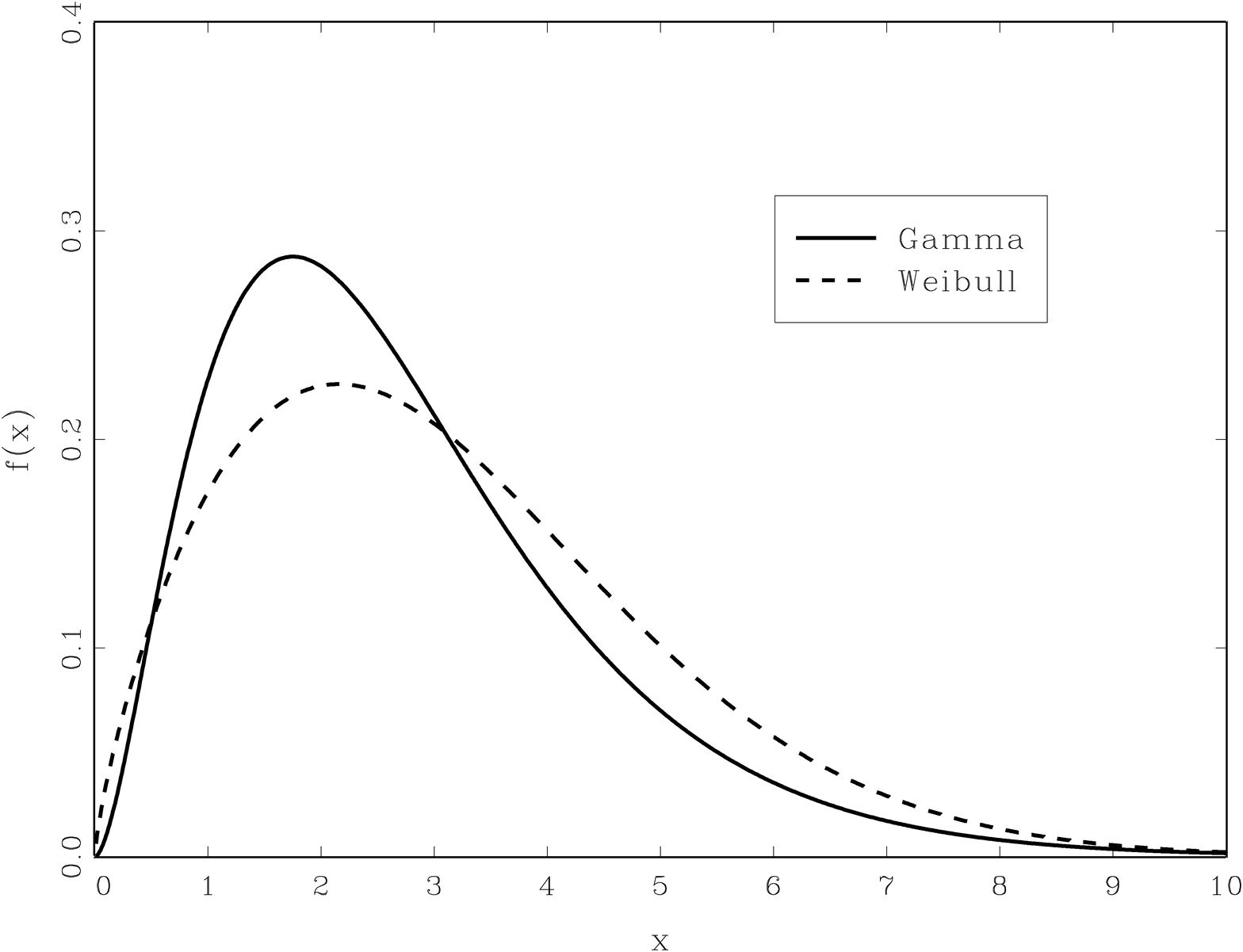}%
}
\bigskip

\textbf{Figure 2:} Density Estimates of School Enrollments\bigskip\bigskip%

{\includegraphics[
height=4.6334in,
width=6.044in
]%
{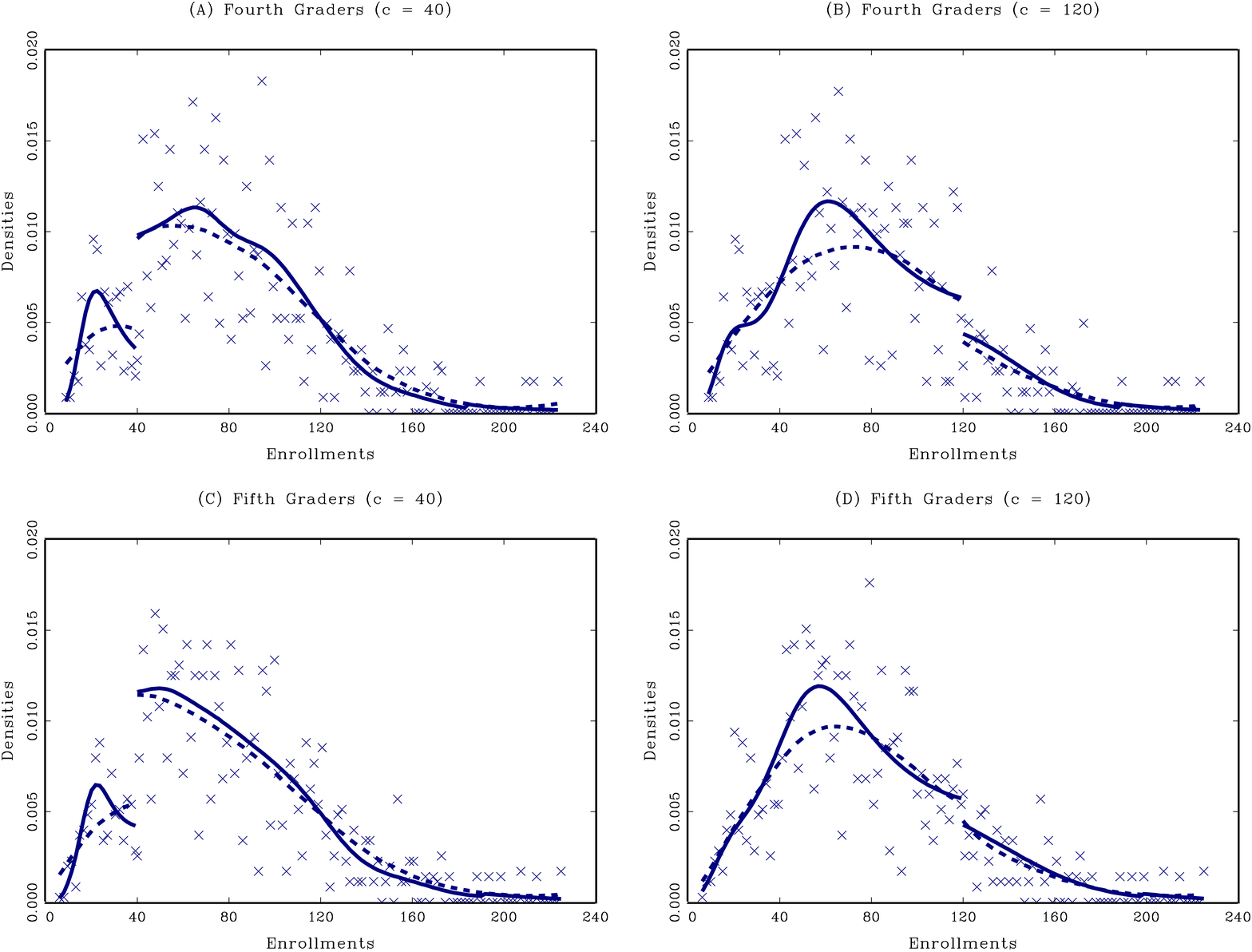}%
}

\end{center}

\paragraph{\textit{Note.}}

In each panel, solid and dashed lines are density estimates via the truncated
gamma kernels and the binned local linear method, respectively. \ The
\textquotedblleft$\times$\textquotedblright\ symbols indicate binned data points.

\begin{center}
\pagebreak

\textbf{Table 1:} Biases, Standard Deviations and RMSEs of Estimators of
$J\left(  c\right)  \bigskip$%

\begin{tabular}
[c]{cccrcccccccc}\hline\hline
&  &  &  &  & \multicolumn{7}{c}{Estimator}\\\cline{6-12}
&  & \multicolumn{3}{r}{} &  &  & \multicolumn{5}{c}{$\tilde{J}\left(
c\right)  $ with $\delta$}\\\cline{8-12}%
Distribution & $c$ & $n$ &  &  & $\hat{J}_{M}\left(  c\right)  $ &  & 0.49 &
& 0.64 &  & 0.81\\\hline
Gamma & 1.7057 & \multicolumn{1}{r}{500} & Bias &  &
\multicolumn{1}{r}{-0.0381} & \multicolumn{1}{r}{} &
\multicolumn{1}{r}{0.0019} & \multicolumn{1}{r}{} & \multicolumn{1}{r}{0.0011}
& \multicolumn{1}{r}{} & \multicolumn{1}{r}{0.0006}\\
\multicolumn{1}{l}{} & (30\%) & \multicolumn{1}{r}{} & StdDev &  &
\multicolumn{1}{r}{0.0461} & \multicolumn{1}{r}{} & \multicolumn{1}{r}{0.0786}
& \multicolumn{1}{r}{} & \multicolumn{1}{r}{0.0812} & \multicolumn{1}{r}{} &
\multicolumn{1}{r}{0.0837}\\
\multicolumn{1}{l}{} &  & \multicolumn{1}{r}{} & RMSE &  &
\multicolumn{1}{r}{0.0598} & \multicolumn{1}{r}{} & \multicolumn{1}{r}{0.0786}
& \multicolumn{1}{r}{} & \multicolumn{1}{r}{0.0812} & \multicolumn{1}{r}{} &
\multicolumn{1}{r}{0.0837}\\\cline{3-12}%
\multicolumn{1}{l}{} &  & \multicolumn{1}{r}{1000} & Bias &  &
\multicolumn{1}{r}{-0.0335} & \multicolumn{1}{r}{} &
\multicolumn{1}{r}{0.0019} & \multicolumn{1}{r}{} & \multicolumn{1}{r}{0.0010}
& \multicolumn{1}{r}{} & \multicolumn{1}{r}{0.0006}\\
&  & \multicolumn{1}{r}{} & StdDev &  & \multicolumn{1}{r}{0.0331} &
\multicolumn{1}{r}{} & \multicolumn{1}{r}{0.0588} & \multicolumn{1}{r}{} &
\multicolumn{1}{r}{0.0607} & \multicolumn{1}{r}{} & \multicolumn{1}{r}{0.0626}%
\\
&  & \multicolumn{1}{r}{} & RMSE &  & \multicolumn{1}{r}{0.0471} &
\multicolumn{1}{r}{} & \multicolumn{1}{r}{0.0588} & \multicolumn{1}{r}{} &
\multicolumn{1}{r}{0.0607} & \multicolumn{1}{r}{} & \multicolumn{1}{r}{0.0626}%
\\\cline{3-12}
&  & \multicolumn{1}{r}{2000} & Bias &  & \multicolumn{1}{r}{-0.0283} &
\multicolumn{1}{r}{} & \multicolumn{1}{r}{0.0006} & \multicolumn{1}{r}{} &
\multicolumn{1}{r}{0.0002} & \multicolumn{1}{r}{} &
\multicolumn{1}{r}{-0.0000}\\
&  & \multicolumn{1}{r}{} & StdDev &  & \multicolumn{1}{r}{0.0250} &
\multicolumn{1}{r}{} & \multicolumn{1}{r}{0.0430} & \multicolumn{1}{r}{} &
\multicolumn{1}{r}{0.0445} & \multicolumn{1}{r}{} & \multicolumn{1}{r}{0.0458}%
\\
&  & \multicolumn{1}{r}{} & RMSE &  & \multicolumn{1}{r}{0.0377} &
\multicolumn{1}{r}{} & \multicolumn{1}{r}{0.0430} & \multicolumn{1}{r}{} &
\multicolumn{1}{r}{0.0445} & \multicolumn{1}{r}{} & \multicolumn{1}{r}{0.0458}%
\\\cline{2-12}
& 2.4248 & \multicolumn{1}{r}{500} & Bias &  & \multicolumn{1}{r}{-0.0407} &
\multicolumn{1}{r}{} & \multicolumn{1}{r}{-0.0011} & \multicolumn{1}{r}{} &
\multicolumn{1}{r}{-0.0011} & \multicolumn{1}{r}{} &
\multicolumn{1}{r}{-0.0011}\\
\multicolumn{1}{l}{} & (Med) & \multicolumn{1}{r}{} & StdDev &  &
\multicolumn{1}{r}{0.0480} & \multicolumn{1}{r}{} & \multicolumn{1}{r}{0.0626}
& \multicolumn{1}{r}{} & \multicolumn{1}{r}{0.0648} & \multicolumn{1}{r}{} &
\multicolumn{1}{r}{0.0668}\\
\multicolumn{1}{l}{} & \multicolumn{1}{l}{} & \multicolumn{1}{r}{} & RMSE &  &
\multicolumn{1}{r}{0.0629} & \multicolumn{1}{r}{} & \multicolumn{1}{r}{0.0626}
& \multicolumn{1}{r}{} & \multicolumn{1}{r}{0.0648} & \multicolumn{1}{r}{} &
\multicolumn{1}{r}{0.0668}\\\cline{3-12}%
\multicolumn{1}{l}{} & \multicolumn{1}{l}{} & \multicolumn{1}{r}{1000} &
Bias &  & \multicolumn{1}{r}{-0.0323} & \multicolumn{1}{r}{} &
\multicolumn{1}{r}{-0.0012} & \multicolumn{1}{r}{} &
\multicolumn{1}{r}{-0.0012} & \multicolumn{1}{r}{} &
\multicolumn{1}{r}{-0.0013}\\
\multicolumn{1}{l}{} & \multicolumn{1}{l}{} & \multicolumn{1}{r}{} & StdDev &
& \multicolumn{1}{r}{0.0353} & \multicolumn{1}{r}{} &
\multicolumn{1}{r}{0.0463} & \multicolumn{1}{r}{} & \multicolumn{1}{r}{0.0479}
& \multicolumn{1}{r}{} & \multicolumn{1}{r}{0.0494}\\
\multicolumn{1}{l}{} & \multicolumn{1}{l}{} & \multicolumn{1}{r}{} & RMSE &  &
\multicolumn{1}{r}{0.0479} & \multicolumn{1}{r}{} & \multicolumn{1}{r}{0.0463}
& \multicolumn{1}{r}{} & \multicolumn{1}{r}{0.0479} & \multicolumn{1}{r}{} &
\multicolumn{1}{r}{0.0494}\\\cline{3-12}%
\multicolumn{1}{l}{} & \multicolumn{1}{l}{} & \multicolumn{1}{r}{2000} &
Bias &  & \multicolumn{1}{r}{-0.0240} & \multicolumn{1}{r}{} &
\multicolumn{1}{r}{-0.0004} & \multicolumn{1}{r}{} &
\multicolumn{1}{r}{-0.0004} & \multicolumn{1}{r}{} &
\multicolumn{1}{r}{-0.0004}\\
\multicolumn{1}{l}{} & \multicolumn{1}{l}{} & \multicolumn{1}{l}{} & StdDev &
& \multicolumn{1}{r}{0.0271} & \multicolumn{1}{r}{} &
\multicolumn{1}{r}{0.0351} & \multicolumn{1}{r}{} & \multicolumn{1}{r}{0.0363}
& \multicolumn{1}{r}{} & \multicolumn{1}{r}{0.0374}\\
\multicolumn{1}{l}{} & \multicolumn{1}{l}{} & \multicolumn{1}{l}{} & RMSE &  &
\multicolumn{1}{r}{0.0362} & \multicolumn{1}{r}{} & \multicolumn{1}{r}{0.0351}
& \multicolumn{1}{r}{} & \multicolumn{1}{r}{0.0363} & \multicolumn{1}{r}{} &
\multicolumn{1}{r}{0.0374}\\\hline
Weibull & 1.9419 & \multicolumn{1}{r}{500} & Bias &  &
\multicolumn{1}{r}{-0.0235} & \multicolumn{1}{r}{} &
\multicolumn{1}{r}{0.0024} & \multicolumn{1}{r}{} & \multicolumn{1}{r}{0.0012}
& \multicolumn{1}{r}{} & \multicolumn{1}{r}{0.0005}\\
& (30\%) & \multicolumn{1}{r}{} & StdDev &  & \multicolumn{1}{r}{0.0416} &
\multicolumn{1}{r}{} & \multicolumn{1}{r}{0.0665} & \multicolumn{1}{r}{} &
\multicolumn{1}{r}{0.0684} & \multicolumn{1}{r}{} & \multicolumn{1}{r}{0.0704}%
\\
\multicolumn{1}{l}{} & \multicolumn{1}{l}{} & \multicolumn{1}{r}{} & RMSE &  &
\multicolumn{1}{r}{0.0478} & \multicolumn{1}{r}{} & \multicolumn{1}{r}{0.0665}
& \multicolumn{1}{r}{} & \multicolumn{1}{r}{0.0684} & \multicolumn{1}{r}{} &
\multicolumn{1}{r}{0.0704}\\\cline{3-12}%
\multicolumn{1}{l}{} & \multicolumn{1}{l}{} & \multicolumn{1}{r}{1000} &
Bias &  & \multicolumn{1}{r}{-0.0187} & \multicolumn{1}{r}{} &
\multicolumn{1}{r}{0.0035} & \multicolumn{1}{r}{} & \multicolumn{1}{r}{0.0013}
& \multicolumn{1}{r}{} & \multicolumn{1}{r}{0.0005}\\
\multicolumn{1}{l}{} & \multicolumn{1}{l}{} & \multicolumn{1}{r}{} & StdDev &
& \multicolumn{1}{r}{0.0302} & \multicolumn{1}{r}{} &
\multicolumn{1}{r}{0.0500} & \multicolumn{1}{r}{} & \multicolumn{1}{r}{0.0509}
& \multicolumn{1}{r}{} & \multicolumn{1}{r}{0.0523}\\
\multicolumn{1}{l}{} & \multicolumn{1}{l}{} & \multicolumn{1}{r}{} & RMSE &  &
\multicolumn{1}{r}{0.0355} & \multicolumn{1}{r}{} & \multicolumn{1}{r}{0.0502}
& \multicolumn{1}{r}{} & \multicolumn{1}{r}{0.0509} & \multicolumn{1}{r}{} &
\multicolumn{1}{r}{0.0523}\\\cline{3-12}%
\multicolumn{1}{l}{} & \multicolumn{1}{l}{} & \multicolumn{1}{r}{2000} &
Bias &  & \multicolumn{1}{r}{-0.0144} & \multicolumn{1}{r}{} &
\multicolumn{1}{r}{0.0017} & \multicolumn{1}{r}{} & \multicolumn{1}{r}{0.0003}
& \multicolumn{1}{r}{} & \multicolumn{1}{r}{0.0001}\\
\multicolumn{1}{l}{} & \multicolumn{1}{l}{} & \multicolumn{1}{r}{} & StdDev &
& \multicolumn{1}{r}{0.0225} & \multicolumn{1}{r}{} &
\multicolumn{1}{r}{0.0367} & \multicolumn{1}{r}{} & \multicolumn{1}{r}{0.0372}
& \multicolumn{1}{r}{} & \multicolumn{1}{r}{0.0383}\\
\multicolumn{1}{l}{} & \multicolumn{1}{l}{} & \multicolumn{1}{r}{} & RMSE &  &
\multicolumn{1}{r}{0.0267} & \multicolumn{1}{r}{} & \multicolumn{1}{r}{0.0367}
& \multicolumn{1}{r}{} & \multicolumn{1}{r}{0.0372} & \multicolumn{1}{r}{} &
\multicolumn{1}{r}{0.0383}\\\cline{2-12}%
\multicolumn{1}{l}{} & 2.8386 & \multicolumn{1}{r}{500} & Bias &  &
\multicolumn{1}{r}{-0.0246} & \multicolumn{1}{r}{} &
\multicolumn{1}{r}{0.0006} & \multicolumn{1}{r}{} & \multicolumn{1}{r}{0.0004}
& \multicolumn{1}{r}{} & \multicolumn{1}{r}{0.0002}\\
\multicolumn{1}{l}{} & (Med) & \multicolumn{1}{r}{} & StdDev &  &
\multicolumn{1}{r}{0.0405} & \multicolumn{1}{r}{} & \multicolumn{1}{r}{0.0534}
& \multicolumn{1}{r}{} & \multicolumn{1}{r}{0.0552} & \multicolumn{1}{r}{} &
\multicolumn{1}{r}{0.0569}\\
\multicolumn{1}{l}{} & \multicolumn{1}{l}{} & \multicolumn{1}{r}{} & RMSE &  &
\multicolumn{1}{r}{0.0474} & \multicolumn{1}{r}{} & \multicolumn{1}{r}{0.0534}
& \multicolumn{1}{r}{} & \multicolumn{1}{r}{0.0552} & \multicolumn{1}{r}{} &
\multicolumn{1}{r}{0.0569}\\\cline{3-12}%
\multicolumn{1}{l}{} & \multicolumn{1}{l}{} & \multicolumn{1}{r}{1000} &
Bias &  & \multicolumn{1}{r}{-0.0195} & \multicolumn{1}{r}{} &
\multicolumn{1}{r}{0.0002} & \multicolumn{1}{r}{} &
\multicolumn{1}{r}{-0.0001} & \multicolumn{1}{r}{} &
\multicolumn{1}{r}{-0.0003}\\
\multicolumn{1}{l}{} & \multicolumn{1}{l}{} & \multicolumn{1}{r}{} & StdDev &
& \multicolumn{1}{r}{0.0290} & \multicolumn{1}{r}{} &
\multicolumn{1}{r}{0.0394} & \multicolumn{1}{r}{} & \multicolumn{1}{r}{0.0408}
& \multicolumn{1}{r}{} & \multicolumn{1}{r}{0.0421}\\
\multicolumn{1}{l}{} & \multicolumn{1}{l}{} & \multicolumn{1}{r}{} & RMSE &  &
\multicolumn{1}{r}{0.0350} & \multicolumn{1}{r}{} & \multicolumn{1}{r}{0.0394}
& \multicolumn{1}{r}{} & \multicolumn{1}{r}{0.0408} & \multicolumn{1}{r}{} &
\multicolumn{1}{r}{0.0421}\\\cline{3-12}%
\multicolumn{1}{l}{} & \multicolumn{1}{l}{} & \multicolumn{1}{r}{2000} &
Bias &  & \multicolumn{1}{r}{-0.0149} & \multicolumn{1}{r}{} &
\multicolumn{1}{r}{0.0007} & \multicolumn{1}{r}{} & \multicolumn{1}{r}{0.0005}
& \multicolumn{1}{r}{} & \multicolumn{1}{r}{0.0004}\\
\multicolumn{1}{l}{} & \multicolumn{1}{l}{} & \multicolumn{1}{r}{} & StdDev &
& \multicolumn{1}{r}{0.0218} & \multicolumn{1}{r}{} &
\multicolumn{1}{r}{0.0299} & \multicolumn{1}{r}{} & \multicolumn{1}{r}{0.0309}
& \multicolumn{1}{r}{} & \multicolumn{1}{r}{0.0319}\\
\multicolumn{1}{l}{} & \multicolumn{1}{l}{} & \multicolumn{1}{r}{} & RMSE &  &
\multicolumn{1}{r}{0.0264} & \multicolumn{1}{r}{} & \multicolumn{1}{r}{0.0299}
& \multicolumn{1}{r}{} & \multicolumn{1}{r}{0.0309} & \multicolumn{1}{r}{} &
\multicolumn{1}{r}{0.0319}\\\hline\hline
\end{tabular}

\pagebreak

\textbf{Table 2:} Finite-Sample Size Properties of Test Statistics for
Discontinuity\bigskip%

\begin{tabular}
[c]{cccrcccccccc}
&  &  &  &  &  &  &  &  &  &  & $\left(  \%\right)  $\\\hline\hline
&  &  & \multicolumn{1}{r}{} &  & \multicolumn{3}{c}{$T_{1}\left(  c\right)  $
with $\delta$} &  & \multicolumn{3}{c}{$T_{2}\left(  c\right)  $ with $\delta
$}\\\cline{6-8}\cline{10-12}%
Distribution & $c$ & $n$ & \multicolumn{1}{r}{Nominal} & \multicolumn{1}{c}{}
& 0.49 & 0.64 & 0.81 &  & 0.49 & 0.64 & 0.81\\\hline
Gamma & 1.7057 & \multicolumn{1}{r}{500} & 5\% &  & \multicolumn{1}{r}{3.2} &
\multicolumn{1}{r}{3.1} & \multicolumn{1}{r}{3.1} & \multicolumn{1}{r}{} &
\multicolumn{1}{r}{4.4} & \multicolumn{1}{r}{4.0} & \multicolumn{1}{r}{3.7}\\
& (30\%) & \multicolumn{1}{r}{} & 10\% &  & \multicolumn{1}{r}{7.5} &
\multicolumn{1}{r}{7.5} & \multicolumn{1}{r}{7.8} & \multicolumn{1}{r}{} &
\multicolumn{1}{r}{8.8} & \multicolumn{1}{r}{8.8} & \multicolumn{1}{r}{8.6}%
\\\cline{3-12}
&  & \multicolumn{1}{r}{1000} & 5\% &  & \multicolumn{1}{r}{3.9} &
\multicolumn{1}{r}{3.9} & \multicolumn{1}{r}{3.9} & \multicolumn{1}{r}{} &
\multicolumn{1}{r}{6.1} & \multicolumn{1}{r}{4.6} & \multicolumn{1}{r}{4.4}\\
&  & \multicolumn{1}{r}{} & 10\% &  & \multicolumn{1}{r}{8.4} &
\multicolumn{1}{r}{8.2} & \multicolumn{1}{r}{8.2} & \multicolumn{1}{r}{} &
\multicolumn{1}{r}{10.7} & \multicolumn{1}{r}{9.2} & \multicolumn{1}{r}{8.9}%
\\\cline{3-12}
&  & \multicolumn{1}{r}{2000} & 5\% &  & \multicolumn{1}{r}{3.5} &
\multicolumn{1}{r}{3.6} & \multicolumn{1}{r}{3.7} & \multicolumn{1}{r}{} &
\multicolumn{1}{r}{4.2} & \multicolumn{1}{r}{3.9} & \multicolumn{1}{r}{3.9}\\
&  & \multicolumn{1}{r}{} & 10\% &  & \multicolumn{1}{r}{8.1} &
\multicolumn{1}{r}{8.2} & \multicolumn{1}{r}{8.4} & \multicolumn{1}{r}{} &
\multicolumn{1}{r}{8.8} & \multicolumn{1}{r}{8.5} & \multicolumn{1}{r}{8.7}%
\\\cline{2-12}
& 2.4248 & \multicolumn{1}{r}{500} & 5\% &  & \multicolumn{1}{r}{3.3} &
\multicolumn{1}{r}{3.6} & \multicolumn{1}{r}{3.6} & \multicolumn{1}{r}{} &
\multicolumn{1}{r}{3.8} & \multicolumn{1}{r}{3.9} & \multicolumn{1}{r}{4.0}\\
\multicolumn{1}{l}{} & (Med) & \multicolumn{1}{r}{} & 10\% &  &
\multicolumn{1}{r}{7.9} & \multicolumn{1}{r}{7.8} & \multicolumn{1}{r}{7.7} &
\multicolumn{1}{r}{} & \multicolumn{1}{r}{8.7} & \multicolumn{1}{r}{8.6} &
\multicolumn{1}{r}{8.5}\\\cline{3-12}%
\multicolumn{1}{l}{} & \multicolumn{1}{l}{} & \multicolumn{1}{r}{1000} & 5\% &
& \multicolumn{1}{r}{3.7} & \multicolumn{1}{r}{3.8} & \multicolumn{1}{r}{3.9}
& \multicolumn{1}{r}{} & \multicolumn{1}{r}{4.1} & \multicolumn{1}{r}{4.2} &
\multicolumn{1}{r}{4.3}\\
\multicolumn{1}{l}{} & \multicolumn{1}{l}{} & \multicolumn{1}{r}{} & 10\% &  &
\multicolumn{1}{r}{8.0} & \multicolumn{1}{r}{8.2} & \multicolumn{1}{r}{8.0} &
\multicolumn{1}{r}{} & \multicolumn{1}{r}{8.6} & \multicolumn{1}{r}{8.6} &
\multicolumn{1}{r}{8.6}\\\cline{3-12}%
\multicolumn{1}{l}{} & \multicolumn{1}{l}{} & \multicolumn{1}{r}{2000} & 5\% &
& \multicolumn{1}{r}{4.7} & \multicolumn{1}{r}{4.7} & \multicolumn{1}{r}{4.8}
& \multicolumn{1}{r}{} & \multicolumn{1}{r}{4.9} & \multicolumn{1}{r}{5.0} &
\multicolumn{1}{r}{5.1}\\
\multicolumn{1}{l}{} & \multicolumn{1}{l}{} & \multicolumn{1}{r}{} & 10\% &  &
\multicolumn{1}{r}{8.8} & \multicolumn{1}{r}{8.9} & \multicolumn{1}{r}{9.0} &
\multicolumn{1}{r}{} & \multicolumn{1}{r}{9.4} & \multicolumn{1}{r}{9.4} &
\multicolumn{1}{r}{9.5}\\\hline
Weibull & 1.9419 & \multicolumn{1}{r}{500} & 5\% &  & \multicolumn{1}{r}{3.2}
& \multicolumn{1}{r}{3.2} & \multicolumn{1}{r}{3.3} & \multicolumn{1}{r}{} &
\multicolumn{1}{r}{6.2} & \multicolumn{1}{r}{4.9} & \multicolumn{1}{r}{4.1}\\
& (30\%) & \multicolumn{1}{r}{} & 10\% &  & \multicolumn{1}{r}{7.7} &
\multicolumn{1}{r}{7.8} & \multicolumn{1}{r}{7.9} & \multicolumn{1}{r}{} &
\multicolumn{1}{r}{10.7} & \multicolumn{1}{r}{9.4} & \multicolumn{1}{r}{9.0}%
\\\cline{3-12}%
\multicolumn{1}{l}{} & \multicolumn{1}{l}{} & \multicolumn{1}{r}{1000} & 5\% &
& \multicolumn{1}{r}{4.0} & \multicolumn{1}{r}{4.2} & \multicolumn{1}{r}{4.2}
& \multicolumn{1}{r}{} & \multicolumn{1}{r}{10.2} & \multicolumn{1}{r}{6.4} &
\multicolumn{1}{r}{5.2}\\
\multicolumn{1}{l}{} & \multicolumn{1}{l}{} & \multicolumn{1}{r}{} & 10\% &  &
\multicolumn{1}{r}{8.2} & \multicolumn{1}{r}{8.3} & \multicolumn{1}{r}{8.4} &
\multicolumn{1}{r}{} & \multicolumn{1}{r}{14.7} & \multicolumn{1}{r}{10.7} &
\multicolumn{1}{r}{9.4}\\\cline{3-12}%
\multicolumn{1}{l}{} & \multicolumn{1}{l}{} & \multicolumn{1}{r}{2000} & 5\% &
& \multicolumn{1}{r}{3.8} & \multicolumn{1}{r}{3.7} & \multicolumn{1}{r}{3.8}
& \multicolumn{1}{r}{} & \multicolumn{1}{r}{7.7} & \multicolumn{1}{r}{4.4} &
\multicolumn{1}{r}{4.0}\\
\multicolumn{1}{l}{} & \multicolumn{1}{l}{} & \multicolumn{1}{r}{} & 10\% &  &
\multicolumn{1}{r}{8.3} & \multicolumn{1}{r}{8.4} & \multicolumn{1}{r}{8.3} &
\multicolumn{1}{r}{} & \multicolumn{1}{r}{12.4} & \multicolumn{1}{r}{9.0} &
\multicolumn{1}{r}{8.5}\\\cline{2-12}%
\multicolumn{1}{l}{} & 2.8386 & \multicolumn{1}{r}{500} & 5\% &  &
\multicolumn{1}{r}{3.6} & \multicolumn{1}{r}{3.6} & \multicolumn{1}{r}{3.5} &
\multicolumn{1}{r}{} & \multicolumn{1}{r}{3.9} & \multicolumn{1}{r}{4.0} &
\multicolumn{1}{r}{3.9}\\
\multicolumn{1}{l}{} & (Med) & \multicolumn{1}{r}{} & 10\% &  &
\multicolumn{1}{r}{7.8} & \multicolumn{1}{r}{7.7} & \multicolumn{1}{r}{7.7} &
\multicolumn{1}{r}{} & \multicolumn{1}{r}{8.5} & \multicolumn{1}{r}{8.5} &
\multicolumn{1}{r}{8.3}\\\cline{3-12}%
\multicolumn{1}{l}{} & \multicolumn{1}{l}{} & \multicolumn{1}{r}{1000} & 5\% &
& \multicolumn{1}{r}{3.7} & \multicolumn{1}{r}{3.8} & \multicolumn{1}{r}{3.8}
& \multicolumn{1}{r}{} & \multicolumn{1}{r}{4.0} & \multicolumn{1}{r}{4.2} &
\multicolumn{1}{r}{4.2}\\
\multicolumn{1}{l}{} & \multicolumn{1}{l}{} & \multicolumn{1}{r}{} & 10\% &  &
\multicolumn{1}{r}{8.1} & \multicolumn{1}{r}{8.1} & \multicolumn{1}{r}{8.2} &
\multicolumn{1}{r}{} & \multicolumn{1}{r}{8.7} & \multicolumn{1}{r}{8.4} &
\multicolumn{1}{r}{8.6}\\\cline{3-12}%
\multicolumn{1}{l}{} & \multicolumn{1}{l}{} & \multicolumn{1}{r}{2000} & 5\% &
& \multicolumn{1}{r}{4.7} & \multicolumn{1}{r}{4.7} & \multicolumn{1}{r}{4.8}
& \multicolumn{1}{r}{} & \multicolumn{1}{r}{4.9} & \multicolumn{1}{r}{5.0} &
\multicolumn{1}{r}{5.0}\\
\multicolumn{1}{l}{} & \multicolumn{1}{l}{} & \multicolumn{1}{r}{} & 10\% &  &
\multicolumn{1}{r}{8.9} & \multicolumn{1}{r}{9.0} & \multicolumn{1}{r}{9.2} &
\multicolumn{1}{r}{} & \multicolumn{1}{r}{9.4} & \multicolumn{1}{r}{9.4} &
\multicolumn{1}{r}{9.6}\\\hline\hline
\end{tabular}

\pagebreak

\textbf{Table 3:} Finite-Sample Power Properties of Test Statistics for
Discontinuity\bigskip%

\begin{tabular}
[c]{ccrrccccccc}%
\multicolumn{4}{l}{(A) Gamma Distribution} &  &  &  &  &  &  & $\left(
\%\right)  $\\\hline\hline
&  &  &  &  & \multicolumn{6}{c}{$d$}\\\cline{6-8}\cline{9-11}%
$c$ & $n$ & \multicolumn{1}{c}{Test} & Nominal &  & 0.00 & 0.02 & 0.04 &
0.06 & 0.08 & 0.10\\\hline
1.7057 & \multicolumn{1}{r}{500} & \multicolumn{1}{c}{$T_{M}\left(  c\right)
$} & 5\% &  & \multicolumn{1}{r}{4.6} & \multicolumn{1}{r}{1.4} &
\multicolumn{1}{r}{2.2} & \multicolumn{1}{r}{9.2} & \multicolumn{1}{r}{28.9} &
\multicolumn{1}{r}{55.1}\\
(30\%) & \multicolumn{1}{r}{} & \multicolumn{1}{c}{} & 10\% &  &
\multicolumn{1}{r}{10.1} & \multicolumn{1}{r}{4.1} & \multicolumn{1}{r}{6.1} &
\multicolumn{1}{r}{17.9} & \multicolumn{1}{r}{43.5} & \multicolumn{1}{r}{68.7}%
\\
& \multicolumn{1}{r}{} & \multicolumn{1}{c}{$T_{1}\left(  c\right)  $} & 5\% &
& \multicolumn{1}{r}{3.1} & \multicolumn{1}{r}{4.5} & \multicolumn{1}{r}{10.7}
& \multicolumn{1}{r}{17.9} & \multicolumn{1}{r}{85.9} &
\multicolumn{1}{r}{98.1}\\
&  & \multicolumn{1}{c}{} & 10\% &  & \multicolumn{1}{r}{7.8} &
\multicolumn{1}{r}{9.4} & \multicolumn{1}{r}{17.0} & \multicolumn{1}{r}{26.8}
& \multicolumn{1}{r}{88.5} & \multicolumn{1}{r}{99.2}\\
&  & \multicolumn{1}{c}{$T_{2}\left(  c\right)  $} & 5\% &  &
\multicolumn{1}{r}{3.7} & \multicolumn{1}{r}{14.2} & \multicolumn{1}{r}{44.7}
& \multicolumn{1}{r}{60.7} & \multicolumn{1}{r}{93.8} &
\multicolumn{1}{r}{98.8}\\
&  &  & 10\% &  & \multicolumn{1}{r}{8.6} & \multicolumn{1}{r}{19.0} &
\multicolumn{1}{r}{53.6} & \multicolumn{1}{r}{63.9} & \multicolumn{1}{r}{97.2}
& \multicolumn{1}{r}{99.7}\\\cline{2-11}
& 1000 & \multicolumn{1}{c}{$T_{M}\left(  c\right)  $} & 5\% &  &
\multicolumn{1}{r}{6.9} & \multicolumn{1}{r}{1.5} & \multicolumn{1}{r}{4.4} &
\multicolumn{1}{r}{22.0} & \multicolumn{1}{r}{58.7} & \multicolumn{1}{r}{87.5}%
\\
& \multicolumn{1}{r}{} & \multicolumn{1}{c}{} & 10\% &  &
\multicolumn{1}{r}{13.6} & \multicolumn{1}{r}{4.1} & \multicolumn{1}{r}{9.5} &
\multicolumn{1}{r}{35.3} & \multicolumn{1}{r}{72.9} & \multicolumn{1}{r}{93.1}%
\\
& \multicolumn{1}{r}{} & \multicolumn{1}{c}{$T_{1}\left(  c\right)  $} & 5\% &
& \multicolumn{1}{r}{3.9} & \multicolumn{1}{r}{6.6} & \multicolumn{1}{r}{12.8}
& \multicolumn{1}{r}{37.1} & \multicolumn{1}{r}{98.7} &
\multicolumn{1}{r}{100.0}\\
&  & \multicolumn{1}{c}{} & 10\% &  & \multicolumn{1}{r}{8.2} &
\multicolumn{1}{r}{12.4} & \multicolumn{1}{r}{21.4} & \multicolumn{1}{r}{46.0}
& \multicolumn{1}{r}{99.0} & \multicolumn{1}{r}{100.0}\\
&  & \multicolumn{1}{c}{$T_{2}\left(  c\right)  $} & 5\% &  &
\multicolumn{1}{r}{4.4} & \multicolumn{1}{r}{13.9} & \multicolumn{1}{r}{50.3}
& \multicolumn{1}{r}{90.8} & \multicolumn{1}{r}{99.5} &
\multicolumn{1}{r}{100.0}\\
&  &  & 10\% &  & \multicolumn{1}{r}{8.9} & \multicolumn{1}{r}{19.2} &
\multicolumn{1}{r}{54.9} & \multicolumn{1}{r}{92.7} & \multicolumn{1}{r}{99.9}
& \multicolumn{1}{r}{100.0}\\\cline{2-11}
& \multicolumn{1}{r}{2000} & \multicolumn{1}{c}{$T_{M}\left(  c\right)  $} &
5\% &  & \multicolumn{1}{r}{9.9} & \multicolumn{1}{r}{1.6} &
\multicolumn{1}{r}{11.7} & \multicolumn{1}{r}{52.1} & \multicolumn{1}{r}{90.1}
& \multicolumn{1}{r}{99.5}\\
&  & \multicolumn{1}{c}{} & 10\% &  & \multicolumn{1}{r}{19.2} &
\multicolumn{1}{r}{4.5} & \multicolumn{1}{r}{21.1} & \multicolumn{1}{r}{66.6}
& \multicolumn{1}{r}{95.6} & \multicolumn{1}{r}{99.8}\\
&  & \multicolumn{1}{c}{$T_{1}\left(  c\right)  $} & 5\% &  &
\multicolumn{1}{r}{3.7} & \multicolumn{1}{r}{8.4} & \multicolumn{1}{r}{36.8} &
\multicolumn{1}{r}{98.8} & \multicolumn{1}{r}{100.0} &
\multicolumn{1}{r}{100.0}\\
&  & \multicolumn{1}{c}{} & 10\% &  & \multicolumn{1}{r}{8.4} &
\multicolumn{1}{r}{15.2} & \multicolumn{1}{r}{44.1} & \multicolumn{1}{r}{99.2}
& \multicolumn{1}{r}{100.0} & \multicolumn{1}{r}{100.0}\\
&  & \multicolumn{1}{c}{$T_{2}\left(  c\right)  $} & 5\% &  &
\multicolumn{1}{r}{3.9} & \multicolumn{1}{r}{25.1} & \multicolumn{1}{r}{90.2}
& \multicolumn{1}{r}{99.5} & \multicolumn{1}{r}{99.9} &
\multicolumn{1}{r}{100.0}\\
&  &  & 10\% &  & \multicolumn{1}{r}{8.7} & \multicolumn{1}{r}{30.4} &
\multicolumn{1}{r}{94.7} & \multicolumn{1}{r}{99.9} & \multicolumn{1}{r}{99.9}
& \multicolumn{1}{r}{100.0}\\\hline
2.4248 & \multicolumn{1}{r}{500} & \multicolumn{1}{c}{$T_{M}\left(  c\right)
$} & 5\% &  & \multicolumn{1}{r}{9.3} & \multicolumn{1}{r}{3.5} &
\multicolumn{1}{r}{2.6} & \multicolumn{1}{r}{4.4} & \multicolumn{1}{r}{9.1} &
\multicolumn{1}{r}{18.7}\\
(Med) & \multicolumn{1}{r}{} & \multicolumn{1}{c}{} & 10\% &  &
\multicolumn{1}{r}{17.4} & \multicolumn{1}{r}{8.8} & \multicolumn{1}{r}{6.0} &
\multicolumn{1}{r}{8.7} & \multicolumn{1}{r}{16.5} & \multicolumn{1}{r}{29.7}%
\\
\multicolumn{1}{l}{} & \multicolumn{1}{r}{} & \multicolumn{1}{c}{$T_{1}\left(
c\right)  $} & 5\% &  & \multicolumn{1}{r}{3.6} & \multicolumn{1}{r}{4.4} &
\multicolumn{1}{r}{7.1} & \multicolumn{1}{r}{12.3} & \multicolumn{1}{r}{20.6}
& \multicolumn{1}{r}{30.1}\\
\multicolumn{1}{l}{} &  & \multicolumn{1}{c}{} & 10\% &  &
\multicolumn{1}{r}{7.7} & \multicolumn{1}{r}{9.1} & \multicolumn{1}{r}{13.4} &
\multicolumn{1}{r}{21.1} & \multicolumn{1}{r}{31.0} & \multicolumn{1}{r}{42.5}%
\\
\multicolumn{1}{l}{} &  & \multicolumn{1}{c}{$T_{2}\left(  c\right)  $} &
5\% &  & \multicolumn{1}{r}{4.0} & \multicolumn{1}{r}{4.8} &
\multicolumn{1}{r}{7.7} & \multicolumn{1}{r}{13.6} & \multicolumn{1}{r}{21.9}
& \multicolumn{1}{r}{32.7}\\
\multicolumn{1}{l}{} &  &  & 10\% &  & \multicolumn{1}{r}{8.5} &
\multicolumn{1}{r}{9.8} & \multicolumn{1}{r}{14.5} & \multicolumn{1}{r}{22.4}
& \multicolumn{1}{r}{32.3} & \multicolumn{1}{r}{44.2}\\\cline{2-11}
& 1000 & \multicolumn{1}{c}{$T_{M}\left(  c\right)  $} & 5\% &  &
\multicolumn{1}{r}{11.5} & \multicolumn{1}{r}{4.0} & \multicolumn{1}{r}{3.5} &
\multicolumn{1}{r}{9.2} & \multicolumn{1}{r}{23.1} & \multicolumn{1}{r}{46.4}%
\\
& \multicolumn{1}{r}{} & \multicolumn{1}{c}{} & 10\% &  &
\multicolumn{1}{r}{20.1} & \multicolumn{1}{r}{8.9} & \multicolumn{1}{r}{7.5} &
\multicolumn{1}{r}{16.0} & \multicolumn{1}{r}{36.1} & \multicolumn{1}{r}{61.7}%
\\
\multicolumn{1}{l}{} & \multicolumn{1}{r}{} & \multicolumn{1}{c}{$T_{1}\left(
c\right)  $} & 5\% &  & \multicolumn{1}{r}{3.9} & \multicolumn{1}{r}{5.0} &
\multicolumn{1}{r}{10.4} & \multicolumn{1}{r}{20.7} & \multicolumn{1}{r}{35.5}
& \multicolumn{1}{r}{53.2}\\
\multicolumn{1}{l}{} &  & \multicolumn{1}{c}{} & 10\% &  &
\multicolumn{1}{r}{8.0} & \multicolumn{1}{r}{10.4} & \multicolumn{1}{r}{18.3}
& \multicolumn{1}{r}{32.1} & \multicolumn{1}{r}{49.0} &
\multicolumn{1}{r}{65.8}\\
\multicolumn{1}{l}{} &  & \multicolumn{1}{c}{$T_{2}\left(  c\right)  $} &
5\% &  & \multicolumn{1}{r}{4.3} & \multicolumn{1}{r}{5.6} &
\multicolumn{1}{r}{11.3} & \multicolumn{1}{r}{22.2} & \multicolumn{1}{r}{36.8}
& \multicolumn{1}{r}{55.0}\\
\multicolumn{1}{l}{} &  &  & 10\% &  & \multicolumn{1}{r}{8.6} &
\multicolumn{1}{r}{11.1} & \multicolumn{1}{r}{19.2} & \multicolumn{1}{r}{33.1}
& \multicolumn{1}{r}{50.3} & \multicolumn{1}{r}{67.0}\\\cline{2-11}
& \multicolumn{1}{r}{2000} & \multicolumn{1}{c}{$T_{M}\left(  c\right)  $} &
5\% &  & \multicolumn{1}{r}{12.0} & \multicolumn{1}{r}{3.6} &
\multicolumn{1}{r}{7.1} & \multicolumn{1}{r}{23.9} & \multicolumn{1}{r}{55.6}
& \multicolumn{1}{r}{83.9}\\
&  & \multicolumn{1}{c}{} & 10\% &  & \multicolumn{1}{r}{20.7} &
\multicolumn{1}{r}{8.0} & \multicolumn{1}{r}{13.8} & \multicolumn{1}{r}{36.0}
& \multicolumn{1}{r}{68.3} & \multicolumn{1}{r}{91.6}\\
\multicolumn{1}{l}{} &  & \multicolumn{1}{c}{$T_{1}\left(  c\right)  $} &
5\% &  & \multicolumn{1}{r}{4.8} & \multicolumn{1}{r}{7.7} &
\multicolumn{1}{r}{18.1} & \multicolumn{1}{r}{37.8} & \multicolumn{1}{r}{60.7}
& \multicolumn{1}{r}{80.2}\\
\multicolumn{1}{l}{} &  & \multicolumn{1}{c}{} & 10\% &  &
\multicolumn{1}{r}{9.0} & \multicolumn{1}{r}{13.9} & \multicolumn{1}{r}{28.3}
& \multicolumn{1}{r}{50.2} & \multicolumn{1}{r}{72.6} &
\multicolumn{1}{r}{87.9}\\
\multicolumn{1}{l}{} &  & \multicolumn{1}{c}{$T_{2}\left(  c\right)  $} &
5\% &  & \multicolumn{1}{r}{5.1} & \multicolumn{1}{r}{8.2} &
\multicolumn{1}{r}{18.9} & \multicolumn{1}{r}{38.8} & \multicolumn{1}{r}{61.7}
& \multicolumn{1}{r}{85.5}\\
\multicolumn{1}{l}{} &  &  & 10\% &  & \multicolumn{1}{r}{9.5} &
\multicolumn{1}{r}{14.6} & \multicolumn{1}{r}{29.2} & \multicolumn{1}{r}{51.4}
& \multicolumn{1}{r}{73.5} & \multicolumn{1}{r}{90.7}\\\hline\hline
\end{tabular}

\pagebreak

\textbf{Table 3} \textit{(Continued)\bigskip}%

\begin{tabular}
[c]{ccrrccccccc}%
\multicolumn{4}{l}{(B) Weibull Distribution} &  &  &  &  &  &  & $\left(
\%\right)  $\\\hline\hline
&  &  &  &  & \multicolumn{6}{c}{$d$}\\\cline{6-8}\cline{9-11}%
$c$ & $n$ & \multicolumn{1}{c}{Test} & Nominal &  & 0.00 & 0.02 & 0.04 &
0.06 & 0.08 & 0.10\\\hline
1.9419 & \multicolumn{1}{r}{500} & \multicolumn{1}{c}{$T_{M}\left(  c\right)
$} & 5\% &  & \multicolumn{1}{r}{3.4} & \multicolumn{1}{r}{1.9} &
\multicolumn{1}{r}{4.0} & \multicolumn{1}{r}{11.8} & \multicolumn{1}{r}{28.5}
& \multicolumn{1}{r}{48.8}\\
(30\%) & \multicolumn{1}{r}{} & \multicolumn{1}{c}{} & 10\% &  &
\multicolumn{1}{r}{7.6} & \multicolumn{1}{r}{5.0} & \multicolumn{1}{r}{9.0} &
\multicolumn{1}{r}{21.2} & \multicolumn{1}{r}{41.2} & \multicolumn{1}{r}{61.8}%
\\
& \multicolumn{1}{r}{} & \multicolumn{1}{c}{$T_{1}\left(  c\right)  $} & 5\% &
& \multicolumn{1}{r}{3.3} & \multicolumn{1}{r}{4.9} & \multicolumn{1}{r}{14.4}
& \multicolumn{1}{r}{19.1} & \multicolumn{1}{r}{85.0} &
\multicolumn{1}{r}{97.1}\\
&  & \multicolumn{1}{c}{} & 10\% &  & \multicolumn{1}{r}{7.9} &
\multicolumn{1}{r}{9.5} & \multicolumn{1}{r}{19.8} & \multicolumn{1}{r}{26.9}
& \multicolumn{1}{r}{88.8} & \multicolumn{1}{r}{98.6}\\
&  & \multicolumn{1}{c}{$T_{2}\left(  c\right)  $} & 5\% &  &
\multicolumn{1}{r}{4.1} & \multicolumn{1}{r}{16.7} & \multicolumn{1}{r}{42.4}
& \multicolumn{1}{r}{57.5} & \multicolumn{1}{r}{90.6} &
\multicolumn{1}{r}{97.8}\\
&  &  & 10\% &  & \multicolumn{1}{r}{9.0} & \multicolumn{1}{r}{21.8} &
\multicolumn{1}{r}{53.0} & \multicolumn{1}{r}{61.0} & \multicolumn{1}{r}{95.5}
& \multicolumn{1}{r}{99.2}\\\cline{2-11}
& 1000 & \multicolumn{1}{c}{$T_{M}\left(  c\right)  $} & 5\% &  &
\multicolumn{1}{r}{4.2} & \multicolumn{1}{r}{2.2} & \multicolumn{1}{r}{7.5} &
\multicolumn{1}{r}{26.3} & \multicolumn{1}{r}{55.0} & \multicolumn{1}{r}{79.3}%
\\
& \multicolumn{1}{r}{} & \multicolumn{1}{c}{} & 10\% &  &
\multicolumn{1}{r}{8.8} & \multicolumn{1}{r}{5.2} & \multicolumn{1}{r}{15.4} &
\multicolumn{1}{r}{39.3} & \multicolumn{1}{r}{68.2} & \multicolumn{1}{r}{87.4}%
\\
& \multicolumn{1}{r}{} & \multicolumn{1}{c}{$T_{1}\left(  c\right)  $} & 5\% &
& \multicolumn{1}{r}{4.2} & \multicolumn{1}{r}{6.5} & \multicolumn{1}{r}{12.9}
& \multicolumn{1}{r}{42.1} & \multicolumn{1}{r}{98.4} &
\multicolumn{1}{r}{99.9}\\
&  & \multicolumn{1}{c}{} & 10\% &  & \multicolumn{1}{r}{8.4} &
\multicolumn{1}{r}{12.4} & \multicolumn{1}{r}{21.2} & \multicolumn{1}{r}{49.1}
& \multicolumn{1}{r}{99.1} & \multicolumn{1}{r}{100.0}\\
&  & \multicolumn{1}{c}{$T_{2}\left(  c\right)  $} & 5\% &  &
\multicolumn{1}{r}{5.2} & \multicolumn{1}{r}{17.5} & \multicolumn{1}{r}{51.0}
& \multicolumn{1}{r}{88.5} & \multicolumn{1}{r}{98.9} &
\multicolumn{1}{r}{99.9}\\
&  &  & 10\% &  & \multicolumn{1}{r}{9.4} & \multicolumn{1}{r}{23.2} &
\multicolumn{1}{r}{56.5} & \multicolumn{1}{r}{91.1} & \multicolumn{1}{r}{99.7}
& \multicolumn{1}{r}{100.0}\\\cline{2-11}
& \multicolumn{1}{r}{2000} & \multicolumn{1}{c}{$T_{M}\left(  c\right)  $} &
5\% &  & \multicolumn{1}{r}{4.5} & \multicolumn{1}{r}{3.1} &
\multicolumn{1}{r}{18.2} & \multicolumn{1}{r}{53.7} & \multicolumn{1}{r}{84.9}
& \multicolumn{1}{r}{97.5}\\
&  & \multicolumn{1}{c}{} & 10\% &  & \multicolumn{1}{r}{9.5} &
\multicolumn{1}{r}{7.1} & \multicolumn{1}{r}{30.1} & \multicolumn{1}{r}{66.6}
& \multicolumn{1}{r}{91.5} & \multicolumn{1}{r}{99.0}\\
&  & \multicolumn{1}{c}{$T_{1}\left(  c\right)  $} & 5\% &  &
\multicolumn{1}{r}{3.8} & \multicolumn{1}{r}{8.3} & \multicolumn{1}{r}{53.7} &
\multicolumn{1}{r}{98.8} & \multicolumn{1}{r}{100.0} &
\multicolumn{1}{r}{100.0}\\
&  & \multicolumn{1}{c}{} & 10\% &  & \multicolumn{1}{r}{8.3} &
\multicolumn{1}{r}{14.8} & \multicolumn{1}{r}{58.0} & \multicolumn{1}{r}{99.5}
& \multicolumn{1}{r}{100.0} & \multicolumn{1}{r}{100.0}\\
&  & \multicolumn{1}{c}{$T_{2}\left(  c\right)  $} & 5\% &  &
\multicolumn{1}{r}{4.0} & \multicolumn{1}{r}{33.2} & \multicolumn{1}{r}{87.4}
& \multicolumn{1}{r}{98.8} & \multicolumn{1}{r}{99.9} &
\multicolumn{1}{r}{100.0}\\
&  &  & 10\% &  & \multicolumn{1}{r}{8.5} & \multicolumn{1}{r}{39.1} &
\multicolumn{1}{r}{93.0} & \multicolumn{1}{r}{99.6} &
\multicolumn{1}{r}{100.0} & \multicolumn{1}{r}{100.0}\\\hline
2.8386 & \multicolumn{1}{r}{500} & \multicolumn{1}{c}{$T_{M}\left(  c\right)
$} & 5\% &  & \multicolumn{1}{r}{4.6} & \multicolumn{1}{r}{2.5} &
\multicolumn{1}{r}{3.0} & \multicolumn{1}{r}{6.0} & \multicolumn{1}{r}{12.6} &
\multicolumn{1}{r}{24.4}\\
(Med) & \multicolumn{1}{r}{} & \multicolumn{1}{c}{} & 10\% &  &
\multicolumn{1}{r}{9.4} & \multicolumn{1}{r}{6.1} & \multicolumn{1}{r}{6.6} &
\multicolumn{1}{r}{11.8} & \multicolumn{1}{r}{22.2} & \multicolumn{1}{r}{36.8}%
\\
\multicolumn{1}{l}{} & \multicolumn{1}{r}{} & \multicolumn{1}{c}{$T_{1}\left(
c\right)  $} & 5\% &  & \multicolumn{1}{r}{3.5} & \multicolumn{1}{r}{4.5} &
\multicolumn{1}{r}{7.2} & \multicolumn{1}{r}{12.3} & \multicolumn{1}{r}{19.9}
& \multicolumn{1}{r}{29.1}\\
\multicolumn{1}{l}{} &  & \multicolumn{1}{c}{} & 10\% &  &
\multicolumn{1}{r}{7.7} & \multicolumn{1}{r}{9.2} & \multicolumn{1}{r}{13.7} &
\multicolumn{1}{r}{21.2} & \multicolumn{1}{r}{30.3} & \multicolumn{1}{r}{41.1}%
\\
\multicolumn{1}{l}{} &  & \multicolumn{1}{c}{$T_{2}\left(  c\right)  $} &
5\% &  & \multicolumn{1}{r}{3.9} & \multicolumn{1}{r}{4.8} &
\multicolumn{1}{r}{7.6} & \multicolumn{1}{r}{14.6} & \multicolumn{1}{r}{22.1}
& \multicolumn{1}{r}{43.6}\\
\multicolumn{1}{l}{} &  &  & 10\% &  & \multicolumn{1}{r}{8.3} &
\multicolumn{1}{r}{9.7} & \multicolumn{1}{r}{14.4} & \multicolumn{1}{r}{22.9}
& \multicolumn{1}{r}{32.2} & \multicolumn{1}{r}{51.9}\\\cline{2-11}
& 1000 & \multicolumn{1}{c}{$T_{M}\left(  c\right)  $} & 5\% &  &
\multicolumn{1}{r}{5.7} & \multicolumn{1}{r}{2.5} & \multicolumn{1}{r}{4.8} &
\multicolumn{1}{r}{12.7} & \multicolumn{1}{r}{31.0} & \multicolumn{1}{r}{55.2}%
\\
& \multicolumn{1}{r}{} & \multicolumn{1}{c}{} & 10\% &  &
\multicolumn{1}{r}{11.2} & \multicolumn{1}{r}{6.2} & \multicolumn{1}{r}{9.1} &
\multicolumn{1}{r}{22.5} & \multicolumn{1}{r}{45.0} & \multicolumn{1}{r}{68.7}%
\\
\multicolumn{1}{l}{} & \multicolumn{1}{r}{} & \multicolumn{1}{c}{$T_{1}\left(
c\right)  $} & 5\% &  & \multicolumn{1}{r}{3.8} & \multicolumn{1}{r}{5.1} &
\multicolumn{1}{r}{10.2} & \multicolumn{1}{r}{20.2} & \multicolumn{1}{r}{34.3}
& \multicolumn{1}{r}{50.9}\\
\multicolumn{1}{l}{} &  & \multicolumn{1}{c}{} & 10\% &  &
\multicolumn{1}{r}{8.2} & \multicolumn{1}{r}{10.7} & \multicolumn{1}{r}{18.1}
& \multicolumn{1}{r}{31.5} & \multicolumn{1}{r}{47.3} &
\multicolumn{1}{r}{63.8}\\
\multicolumn{1}{l}{} &  & \multicolumn{1}{c}{$T_{2}\left(  c\right)  $} &
5\% &  & \multicolumn{1}{r}{4.2} & \multicolumn{1}{r}{5.7} &
\multicolumn{1}{r}{11.3} & \multicolumn{1}{r}{21.3} & \multicolumn{1}{r}{36.7}
& \multicolumn{1}{r}{61.7}\\
\multicolumn{1}{l}{} &  &  & 10\% &  & \multicolumn{1}{r}{8.6} &
\multicolumn{1}{r}{11.3} & \multicolumn{1}{r}{19.0} & \multicolumn{1}{r}{32.5}
& \multicolumn{1}{r}{48.9} & \multicolumn{1}{r}{70.6}\\\cline{2-11}
& \multicolumn{1}{r}{2000} & \multicolumn{1}{c}{$T_{M}\left(  c\right)  $} &
5\% &  & \multicolumn{1}{r}{6.6} & \multicolumn{1}{r}{3.1} &
\multicolumn{1}{r}{10.0} & \multicolumn{1}{r}{31.8} & \multicolumn{1}{r}{63.9}
& \multicolumn{1}{r}{87.6}\\
&  & \multicolumn{1}{c}{} & 10\% &  & \multicolumn{1}{r}{12.7} &
\multicolumn{1}{r}{6.7} & \multicolumn{1}{r}{17.7} & \multicolumn{1}{r}{46.0}
& \multicolumn{1}{r}{76.1} & \multicolumn{1}{r}{93.3}\\
\multicolumn{1}{l}{} &  & \multicolumn{1}{c}{$T_{1}\left(  c\right)  $} &
5\% &  & \multicolumn{1}{r}{4.8} & \multicolumn{1}{r}{7.8} &
\multicolumn{1}{r}{18.0} & \multicolumn{1}{r}{36.2} & \multicolumn{1}{r}{58.6}
& \multicolumn{1}{r}{79.8}\\
\multicolumn{1}{l}{} &  & \multicolumn{1}{c}{} & 10\% &  &
\multicolumn{1}{r}{9.2} & \multicolumn{1}{r}{14.1} & \multicolumn{1}{r}{28.1}
& \multicolumn{1}{r}{49.0} & \multicolumn{1}{r}{70.8} &
\multicolumn{1}{r}{87.1}\\
\multicolumn{1}{l}{} &  & \multicolumn{1}{c}{$T_{2}\left(  c\right)  $} &
5\% &  & \multicolumn{1}{r}{5.0} & \multicolumn{1}{r}{8.1} &
\multicolumn{1}{r}{18.6} & \multicolumn{1}{r}{37.8} & \multicolumn{1}{r}{64.7}
& \multicolumn{1}{r}{99.2}\\
\multicolumn{1}{l}{} &  &  & 10\% &  & \multicolumn{1}{r}{9.6} &
\multicolumn{1}{r}{14.6} & \multicolumn{1}{r}{28.7} & \multicolumn{1}{r}{49.9}
& \multicolumn{1}{r}{74.6} & \multicolumn{1}{r}{99.5}\\\hline\hline
\end{tabular}

\end{center}

\paragraph{\textit{Note.}}

The value of $\delta$ for each of $T_{1}\left(  c\right)  $ and $T_{2}\left(
c\right)  $ is set equal to 0.81.

\pagebreak

\begin{center}
\textbf{Table 4:} Estimation and Testing for the Discontinuity of Densities of
School Enrollments\bigskip%

\begin{tabular}
[c]{ccccccccccc}\hline\hline
&  & \multicolumn{4}{c}{Binned Local Linear Method} &  &
\multicolumn{4}{c}{Truncated Kernel Method}\\\cline{3-6}\cline{8-11}%
$n$ & $c$ & $\hat{f}_{-}^{M}\left(  c\right)  $ & $\hat{f}_{+}^{M}\left(
c\right)  $ & $\hat{J}_{M}\left(  c\right)  $ & $T_{M}\left(  c\right)  $ &  &
$\tilde{f}_{-}\left(  c\right)  $ & $\tilde{f}_{+}\left(  c\right)  $ &
$\tilde{J}\left(  c\right)  $ & $T_{2}\left(  c\right)  $\\\hline
&  & \multicolumn{9}{c}{\textit{(a) Fourth Graders:}}\\
\multicolumn{1}{r}{2059} & \multicolumn{1}{r}{40} & \multicolumn{1}{r}{0.0046}
& \multicolumn{1}{r}{0.0096} & \multicolumn{1}{r}{0.0050} &
\multicolumn{1}{r}{\textbf{5.61}} & \multicolumn{1}{r}{} &
\multicolumn{1}{r}{0.0034} & \multicolumn{1}{r}{0.0098} &
\multicolumn{1}{r}{0.0064} & \multicolumn{1}{r}{\textbf{5.76}}\\
\multicolumn{1}{r}{} & \multicolumn{1}{r}{80} & \multicolumn{1}{r}{0.0103} &
\multicolumn{1}{r}{0.0097} & \multicolumn{1}{r}{-0.0006} &
\multicolumn{1}{r}{-0.62} & \multicolumn{1}{r}{} & \multicolumn{1}{r}{0.0086}
& \multicolumn{1}{r}{0.0090} & \multicolumn{1}{r}{0.0003} &
\multicolumn{1}{r}{0.24}\\
\multicolumn{1}{r}{} & \multicolumn{1}{r}{120} & \multicolumn{1}{r}{0.0061} &
\multicolumn{1}{r}{0.0039} & \multicolumn{1}{r}{-0.0022} &
\multicolumn{1}{r}{\textbf{-3.35}} & \multicolumn{1}{r}{} &
\multicolumn{1}{r}{0.0063} & \multicolumn{1}{r}{0.0044} &
\multicolumn{1}{r}{-0.0020} & \multicolumn{1}{r}{\textbf{-3.55}}\\
\multicolumn{1}{r}{} & \multicolumn{1}{r}{160} & \multicolumn{1}{r}{0.0011} &
\multicolumn{1}{r}{0.0009} & \multicolumn{1}{r}{-0.0003} &
\multicolumn{1}{r}{-0.84} & \multicolumn{1}{r}{} & \multicolumn{1}{r}{0.0013}
& \multicolumn{1}{r}{0.0005} & \multicolumn{1}{r}{-0.0008} &
\multicolumn{1}{r}{\textbf{-2.88}}\\\hline
\multicolumn{1}{r}{} & \multicolumn{1}{r}{} & \multicolumn{9}{c}{\textit{(b)
Fifth Graders:}}\\
\multicolumn{1}{r}{2029} & \multicolumn{1}{r}{40} & \multicolumn{1}{r}{0.0055}
& \multicolumn{1}{r}{0.0114} & \multicolumn{1}{r}{0.0059} &
\multicolumn{1}{r}{\textbf{6.29}} & \multicolumn{1}{r}{} &
\multicolumn{1}{r}{0.0042} & \multicolumn{1}{r}{0.0116} &
\multicolumn{1}{r}{0.0074} & \multicolumn{1}{r}{\textbf{6.28}}\\
\multicolumn{1}{r}{} & \multicolumn{1}{r}{80} & \multicolumn{1}{r}{0.0107} &
\multicolumn{1}{r}{0.0098} & \multicolumn{1}{r}{-0.0009} &
\multicolumn{1}{r}{-0.98} & \multicolumn{1}{r}{} & \multicolumn{1}{r}{0.0087}
& \multicolumn{1}{r}{0.0103} & \multicolumn{1}{r}{0.0017} &
\multicolumn{1}{r}{1.25}\\
\multicolumn{1}{r}{} & \multicolumn{1}{r}{120} & \multicolumn{1}{r}{0.0054} &
\multicolumn{1}{r}{0.0045} & \multicolumn{1}{r}{-0.0009} &
\multicolumn{1}{r}{-1.20} & \multicolumn{1}{r}{} & \multicolumn{1}{r}{0.0057}
& \multicolumn{1}{r}{0.0043} & \multicolumn{1}{r}{-0.0014} &
\multicolumn{1}{r}{\textbf{-2.84}}\\
\multicolumn{1}{r}{} & \multicolumn{1}{r}{160} & \multicolumn{1}{r}{0.0014} &
\multicolumn{1}{r}{0.0011} & \multicolumn{1}{r}{-0.0003} &
\multicolumn{1}{r}{-0.80} & \multicolumn{1}{r}{} & \multicolumn{1}{r}{0.0014}
& \multicolumn{1}{r}{0.0010} & \multicolumn{1}{r}{-0.0004} &
\multicolumn{1}{r}{-1.28}\\\hline\hline
\end{tabular}

\end{center}

\paragraph{\textit{Note.}}

The value of $\delta$ for $T_{2}\left(  c\right)  $ is set equal to 0.81.
\ Values of test statistics in bold faces indicate significance at the 5\% level.

\end{document}